\documentclass[11pt]{amsart}%
\usepackage{amsmath}
\usepackage{graphicx}%
\usepackage{amsfonts}%
\usepackage{amssymb}
\newtheorem{theorem}{Theorem}

\newtheorem{corollary}[theorem]{Corollary}

\newtheorem{lemma}[theorem]{Lemma}

\newtheorem{proposition}[theorem]{Proposition}

\begin{document}
\title{The Bourgain $\ell^{1}$-index of mixed Tsirelson space}
\author{Denny H. Leung}
\address{Department of Mathematics \\
National University of Singapore \\
Singapore 117543}
\email{matlhh@nus.edu.sg}
\author{Wee-Kee Tang}
\address{Mathematics and Mathematics Education\\
National Institute of Education\\
Nanyang Technological University\\
1 Nanyang Walk, Singapore 637616}
\email{wktang@nie.edu.sg}

\begin{abstract}
Suppose that $(\mathcal{F}_{n})_{n=0}^{\infty}$ is a sequence of regular
families of finite subsets of $\mathbb{N}$ such that $\mathcal{F}_{0}$
contains all singletons, and $(\theta_{n})_{n=1}^{\infty}$ is a nonincreasing
null sequence in $(0,1)$. The mixed Tsirelson space $T(\mathcal{F}_{0}%
,(\theta_{n},\mathcal{F}_{n})_{n=1}^{\infty})$ is the completion of $c_{00}$
with respect to the implicitly defined norm
\[
\left\|  x\right\|  =\max\left\{  \left\|  x\right\|  _{\mathcal{F}_{0}}%
,\sup\limits_{n\in\mathbb{N}}\sup\theta_{n}\sum_{i=1}^{k}\left\|
E_{i}x\right\|  \right\}  ,
\]
where $\Vert x\Vert_{\mathcal{F}_{0}}=\sup_{F\in\mathcal{F}}\Vert
Fx\Vert_{\ell^{1}}$ and the last supremum is taken over all sequences
$(E_{i})_{i=1}^{k}$ in $[\mathbb{N}]^{<\infty}$ such that $\max E_{i}<\min
E_{i+1}$ and $\{\min E_{i}:1\leq i\leq k\}\in\mathcal{F}_{n}$. In this paper,
we compute the Bourgain $\ell^{1}$-index of the space $T(\mathcal{F}%
_{0},(\theta_{n},\mathcal{F}_{n})_{n=1}^{\infty})$. As a consequence, it is
shown that if $\eta$ is a countable ordinal not of the form $\omega^{\xi}$ for
some limit ordinal $\xi$, then there is a Banach space whose $\ell^{1}$-index
is $\omega^{\eta}$.

\end{abstract}
\maketitle


\section{Introduction}

Endow the power set of $\mathbb{N}$, identified with $2^{\mathbb{N}}$, with
the product topology. Denote by $[\mathbb{N}]^{<\infty}$ the subspace
consisting of all finite subsets of $\mathbb{N}$. A family $\mathcal{F}%
\subseteq\lbrack\mathbb{N}]^{<\infty}$ is said to be \emph{hereditary} if
$G\subseteq F\in\mathcal{F}$ implies $G\in\mathcal{F}$. It is \emph{spreading}
if whenever $F=\{n_{1},\dots,n_{k}\}\in\mathcal{F}$, $n_{1}<\dots<n_{k}$, and
$m_{1}<\dots<m_{k}$ satisfy $m_{i}\geq n_{i}$, $1\leq i\leq k$, then
$\{m_{1},\dots,m_{k}\}\in\mathcal{F}$. In this case, we also say that
$\{m_{1},\dots,m_{k}\}$ is a \emph{spreading} of $F$. A \emph{regular} family
is one that is hereditary, spreading and compact (as a subset of the
topological space $[\mathbb{N}]^{<\infty}$). Let $c_{00}$ be the vector space
of all finitely supported real sequences and let $(e_{k})$ be the standard
unit vector basis of $c_{00}$. If $\mathcal{F}$ is regular, define the
seminorm $\Vert\cdot\Vert_{\mathcal{F}}$ on $c_{00}$
by $\Vert\sum a_{k}e_{k}\Vert_{\mathcal{F}}=\sup_{F\in\mathcal{F}}\sum_{k\in
F}|a_{k}|$. For $E\in\lbrack\mathbb{N}]^{<\infty}$ and $x = \sum a_{k}e_{k}\in
c_{00}$, let $Ex = \sum_{k\in E}a_{k}e_{k}\in c_{00}$.
Given a sequence of regular families $(\mathcal{F}_{n})_{n=0}^{\infty}$ such
that $\mathcal{F}_{0}$ contains all singleton subsets of $\mathbb{N}$, and a
nonincreasing null sequence $\left(  \theta_{n}\right)  _{n=1}^{\infty}$ in
$\left(  0,1\right)  $, the \emph{mixed Tsirelson space} $T\left(
\mathcal{F}_{0},\left(  \theta_{n},\mathcal{F}_{n}\right)  _{n=1}^{\infty
}\right)  $ is the completion of $c_{00}$ under the implicitly defined norm%
\begin{equation}
\left\|  x\right\|  =\max\left\{  \left\|  x\right\|  _{\mathcal{F}_{0}}%
,\sup\limits_{n\in\mathbb{N}}\sup\theta_{n}\sum_{i=1}^{k}\left\|
E_{i}x\right\|  \right\}  , \label{norm}%
\end{equation}
where the last supremum is taken over all sequences $(E_{i})_{i=1}^{k}$ in
$[\mathbb{N}]^{<\infty}$ such that $\max E_{i}<\min E_{i+1}$ and $\{\min
E_{i}:1\leq i\leq k\}\in\mathcal{F}_{n}$. The main aim of the present paper is
the computation of the $\ell^{1}$-index $I_{b}\left(  T\left(  \mathcal{F}%
_{0},\left(  \theta_{n},\mathcal{F}_{n}\right)  _{n=1}^{\infty}\right)
\right)  $ (defined below) in terms of the sequences $(\mathcal{F}_{n}%
)_{n=0}^{\infty}$ and $\left(  \theta_{n}\right)  _{n=1}^{\infty}.$ It follows
from our work (see Corollary \ref{solution} below) that if $\eta$ is a
countable ordinal not of the form $\omega^{\xi}$ for some limit ordinal $\xi$,
then there is a Banach space whose $\ell^{1}$-index is $\omega^{\eta}$. This
answers Question 1 in \cite{JO}.

Our starting point is a comparison of normalized block basic sequences in
$T(\mathcal{F}_{0},(\theta_{n},\mathcal{F}_{n})_{n=1}^{\infty})$ with
subsequences of the unit vector basis in related mixed Tsirelson spaces
(Proposition \ref{comparison}). In particular, we obtain in Corollary
\ref{equivalence} that every normalized block basic sequence in a mixed
Tsirelson space $T(\mathcal{F}_{0},(\theta_{n},\mathcal{F}_{n})_{n=1}^{\ell})$
defined by finitely many families is equivalent to a subsequence of the unit
vector basis in the same space. This result was proved for the Figiel-Johnson
Tsirelson space in \cite{CJT} and for certain generalized Tsirelson spaces in
\cite{Be}. Our approach may be considered as a descendant of that in \cite{Be}.

In \S3, the comparison result is used to obtain bounds on the $\ell^{1}%
$-index. In \S4, we introduce a method of constructing $\ell^{1}$-trees of
large index. This is a two-step method whereby many $\ell^{1}(n)$-block basic
sequences are first constructed (Lemma \ref{LB3}) and these are then condensed
into $\ell^{1}$-trees by a compactness argument (Lemma \ref{P20}).

If $M$ is an infinite subset of $\mathbb{N}$, denote the set of all finite,
respectively infinite, subsets of $M$ by $[M]^{<\infty},$ respectively
$\left[  M\right]  $. If $E$ and $F$ are finite subsets of $\mathbb{N}$, we
write $E<F$, respectively $E\leq F$, to mean $\max E<\min F$, respectively
$\max E\leq\min F$ ($\max\emptyset=0$ and $\min\emptyset=\infty$). We
abbreviate $\{n\}<E$ and $\{n\}\leq E$ to $n<E$ and $n\leq E$ respectively.
Given $\mbox{$\mathcal{F}$}\subseteq\lbrack\mathbb{N}]^{<\infty}$, a sequence
of finite subsets $\{E_{1},\dots,E_{n}\}$ of $\mathbb{N}$ is said to be
$\mbox{$\mathcal{F}$}$-\emph{admissible} if $E_{1}<\dots<E_{n}$ and $\{\min
E_{1},\dots,\min E_{n}\}\in\mbox{$\mathcal{F}$}$. If $\mbox{$\mathcal{M}$}$
and $\mbox{$\mathcal{N}$}$ are regular subsets of $[\mathbb{N}]^{<\infty}$, we
let
\[
\mbox{$\mathcal{M}$}[\mbox{$\mathcal{N}$}]=\{\cup_{i=1}^{k}F_{i}:F_{i}\in
\mbox{$\mathcal{N}$}\text{ for all $i$ and }\{F_{1},\dots,F_{k}\}\text{ is
$\mbox{$\mathcal{M}$}$-admissible}\}.
\]
Given a sequence of regular families $(\mathcal{M}_{i})$, we define
inductively $[\mathcal{M}_{1},\mathcal{M}_{2}]=\mathcal{M}_{1}[\mathcal{M}%
_{2}]$ and $[\mathcal{M}_{1},\dots,\mathcal{M}_{i+1}]=[\mathcal{M}_{1}%
,\dots,\mathcal{M}_{i}][\mathcal{M}_{i+1}]$. Also, let
\[
\left(  \mathcal{M}_{1},\dots,\mathcal{M}_{k}\right)  =\left\{  \cup_{i=1}%
^{k}M_{i}:M_{i}\in\mathcal{M}_{i},M_{1}<\dots<M_{k}\right\}  .
\]
We abbreviate the $k$-fold construction $(\mathcal{M},\dots,\mathcal{M})$ as
$(\mathcal{M})^{k}$. Of primary importance are the Schreier classes as defined
in \cite{AA}. We will need a slightly extended version of such classes.
Suppose that $g:\mathbb{N}\rightarrow\mathbb{N}$ is a function increasing to
$\infty$. Let $\mbox{$\mathcal{S}$}_{0}^{g}=\{\{n\}:n\in\mathbb{N}%
\}\cup\{\emptyset\}$ and $\mbox{$\mathcal{S}$}_{1}^{g}=\{F\subseteq
\mathbb{N}:|F|\leq g(\min F)\}$. Here $|F|$ denotes the cardinality of $F$.
The higher Schreier classes are defined inductively as follows. $\mbox
{$\mathcal{S}$}_{\alpha+1}^{g}=\mbox{$\mathcal{S}$}_{1}^{g}[\mbox{$\mathcal
{S}$}_{\alpha}^{g}]$ for all $\alpha<\omega_{1}$. If $\alpha$ is a countable
limit ordinal, choose a sequence $(\alpha_{n})$ strictly increasing to
$\alpha$ and set
\[
\mbox{$\mathcal{S}$}_{\alpha}^{g}=\{F:F\in\mbox{$\mathcal{S}$}_{\alpha_{n}%
}^{g}\text{ for some $n\leq g(|F|)$}\}.
\]
If $g$ is the identity function, then we obtain the usual Schreier classes,
and we abbreviate $\mathcal{S}_{\alpha}^{g}$ to $\mathcal{S}_{\alpha}$. It is
clear that $\mbox{$\mathcal{S}$}_{\alpha}^{g}$ is a regular family for all
$\alpha<\omega_{1}$. If $M=(m_{1},m_{2},\dots)$ is a subsequence of
$\mathbb{N}$, let $\mbox{$\mathcal{S}$}_{\alpha}(M)=\{\{m_{i}:i\in
F\}:F\in\mbox{$\mathcal{S}$}_{\alpha}\}$. Since $\mbox{$\mathcal{S}$}_{\alpha
}$ is spreading, $\mbox{$\mathcal{S}$}_{\alpha}(M)\subseteq\mbox{$\mathcal{S}%
$}_{\alpha}$.

The norm in a mixed Tsirelson space can be computed in terms of trees
(\cite{Be}, \cite{OT}). A \emph{tree} in $\left[  \mathbb{N}\right]
^{<\infty}$ is a finite collection of elements $\left(  E_{i}^{m}\right)  ,$
$0\leq m\leq r,$ $1\leq i\leq k\left(  m\right)  ,$ in $\left[  \mathbb{N}%
\right]  ^{<\infty}$ so that for each $m,$ $E_{1}^{m}<E_{2}^{m}<\dots
<E_{k\left(  m\right)  }^{m},$ and that every $E_{i}^{m+1}$ is a subset of
some $E_{j}^{m}.$ The elements $E_{i}^{m}$ are called \emph{nodes} of the
tree. Any node $E_{i}^{m}$ is said to be of \emph{level} $m.$ Nodes at level
$0$ are called \textit{roots. }If $E_{i}^{n}\subseteq E_{j}^{m}$ and $n>m,$ we
say that $E_{i}^{n}$ is a \emph{descendant }of $E_{j}^{m}$ and $E_{j}^{m}$ is
an \emph{ancestor} of $E_{i}^{n}.$ If, in the above notation, $n=m+1,$ then
$E_{i}^{n}$ is said to be an \emph{immediate successor} of $E_{j}^{m},$ and
$E_{j}^{m}$ the \emph{immediate predecessor }of $E_{i}^{n}.$ Nodes with no
descendants are called \emph{terminal nodes} or \emph{leaves} of the tree.
Given a node $E$ in a tree $\mathcal{T},$ denote by $\mathcal{T}_{E}$ the
subtree consisting of the node $E$ together with all its descendants. A tree
$\left(  E_{i}^{m}\right)  ,$ $0\leq m<r,$ $1\leq i\leq k\left(  m\right)  ,$
is $\left(  \mathcal{F}_{n}\right)  $-admissible if $k\left(  0\right)  =1$
and for every $m$ and $i,$ the collection $\left(  E_{j}^{m+1}\right)  $ of
all immediate successors of $E_{i}^{m}$ is an $\mathcal{F}_{n}$-admissible
collection for some $n\in\mathbb{N}.$ Given an $\left(  \mathcal{F}%
_{n}\right)  $-admissible tree $\left(  E_{i}^{m}\right)  ,$ we define the
\emph{history} of the individual nodes inductively as follows. Let $h\left(
E_{1}^{0}\right)  =\left(  0\right)  .$ If $h\left(  E_{i}^{m}\right)  $ has
been defined and the collection $\left(  E_{j}^{m+1}\right)  $ of all
immediate successors of $E_{i}^{m}$ forms an $\mathcal{F}_{n}$-admissible
collection, then define $h\left(  E_{j}^{m+1}\right)  $ to be the $\left(
m+2\right)  $- tuple $\left(  h\left(  E_{i}^{m}\right)  ,n\right)  $ and let
$n\left(  E_{j}^{m+1}\right)  =n$ for each immediate successor \ $E_{j}^{m+1}$
of $E_{i}^{m}.$ Finally, assign $\left(  \left(  \theta_{n}\right)
\text{-compatible}\right)  $ \emph{tags} to the nodes by defining $t\left(
E_{i}^{m}\right)  =\prod_{j=0}^{m}\theta_{n_{j}}$ if $h\left(  E_{i}%
^{m}\right)  =\left(  n_{0},n_{1},\dots,n_{m}\right)  $ $\left(  \theta
_{0}=1\right)  .$ If $x\in c_{00}$ and $\mathcal{T}$ is an $\left(
\mathcal{F}_{n}\right)  $-admissible tree, let $\mathcal{T}x=\sum t\left(
E\right)  \left\|  Ex\right\|  ,$ where the sum is taken over all leaves in
$\mathcal{T}.$ It is easily observed that $\left\|  x\right\|  =\max\left\{
\mathcal{T}x:\mathcal{T}\text{ is an }\left(  \mathcal{F}_{n}\right)
\text{-admissible tree}\right\}  .$ An $\left(  \mathcal{F}_{n}\right)
$-admissible tree is said to be $\emph{complete}$ (for a particular $x\in
c_{00}$) if $\left\|  Ex\right\|  =\left\|  Ex\right\|  _{\mathcal{F}_{0}}$
for every leaf $E$ in $\mathcal{T}.$ Clearly, for every $x\in c_{00},$ there
is a complete tree $\mathcal{T}$ such that $\left\|  x\right\|  =\mathcal{T}%
x.$ Let us observe that if we define $\|x\|$ to be $\sup\sum
t(E)\|Ex\|_{\mathcal{F}_{0}}$, where the $\sup$ is taken over all
$(\mathcal{F}_{n})$-admissible trees $\mathcal{T}$ and the sum is taken over
all leaves $E$ in $\mathcal{T}$, then the resulting norm satisfies the
implicit equation (\ref{norm}).

\begin{proposition}
\label{P0}Let $T\left(  \mathcal{F}_{0},\left(  \theta_{n},\mathcal{F}%
_{n}\right)  _{n=1}^{\infty}\right)  $ be as above. Choose a strictly
increasing sequence of integers $\left(  m_{k}\right)  _{k=0}^{\infty}$ such
that $m_{0}=0$ and $\theta_{m_{k+1}}\leq\frac{1}{2}\theta_{m_{k}}$ for all
$k\in\mathbb{N}.$ If $m_{k-1}<n\leq m_{k},$ let $\mathcal{G}_{n}=\left\{
F\in\mathcal{F}_{n}:k\leq F\right\}  \cup\mathcal{S}_{0}.$ Then $T\left(
\mathcal{F}_{0},\left(  \theta_{n},\mathcal{F}_{n}\right)  _{n=1}^{\infty
}\right)  $ is isomorphic to $T\left(  \mathcal{F}_{0},\left(  \theta
_{n},\mathcal{G}_{n}\right)  _{n=1}^{\infty}\right)  $ via the formal identity.
\end{proposition}

\begin{proof}
Denote the norms on $T\left(  \mathcal{F}_{0},\left(  \theta_{n}%
,\mathcal{F}_{n}\right)  _{n=1}^{\infty}\right)  $ and $T\left(
\mathcal{F}_{0},\left(  \theta_{n},\mathcal{G}_{n}\right)  _{n=1}^{\infty
}\right)  $ by $\left\|  \cdot\right\|  $ and $\left|  \left|  \left|
\cdot\right|  \right|  \right|  $ respectively. Clearly, $\left|  \left|
\left|  x\right|  \right|  \right|  \leq\left\|  x\right\|  $ for all $x\in
c_{00}.$ Given a fixed element $x\in c_{00},$ let $\mathcal{T}^{\mathcal{F}}$
denote a complete $\left(  \mathcal{F}_{n}\right)  $-admissible tree such that
$\left\|  x\right\|  =\mathcal{T}^{\mathcal{F}}x.$ If $F$ is a node of
$\mathcal{T}^{\mathcal{F}}$ other than the root, let $G_{F}=F\cap\lbrack
k,\infty),$ where $k$ is the unique integer such that $m_{k-1}<\max
\{n_{1},\dots,n_{r}\}\leq m_{k},$ $h\left(  F\right)  =\left(  0,n_{1}%
,\dots,n_{r}\right)  .$ If $F$ is the root of $\mathcal{T}^{\mathcal{F}},$ let
$G_{F}=F.$ Then $\mathcal{T}^{\mathcal{G}}=\left\{  G_{F}:F\in\mathcal{T}%
^{\mathcal{F}}\right\}  $ is a $\left(  \mathcal{G}_{n}\right)  $-admissible
tree. For any $r\in\mathbb{N},$ let $\mathcal{L}_{r}$ be the set of level $r$
leaves in $\mathcal{T}^{\mathcal{F}}.$ Arrange the elements in $\mathcal{L}%
_{r}$ from left to right as $F_{1}<F_{2}<\dots<F_{\ell}.$ If $1\leq j\leq
\ell,$ write $h\left(  F_{j}\right)  =\left(  0,n_{j{1}},,\dots,n_{j{r}%
}\right)  $ and determine $k_{j}$ such that $m_{k_{j}-1}<\max\{n_{j{1}}%
,\dots,n_{j{r}}\}\leq m_{k_{j}}.$ If $k_{j}\leq j,$ then $k_{j}\leq j\leq
F_{j}.$ Thus $G_{F_{j}}=F_{j}\cap\lbrack k_{j},\infty)=F_{j}.$ Otherwise,
$j<k_{j},$ and hence
\[
t\left(  F_{j}\right)  \left\|  F_{j}x\right\|  _{\mathcal{F}_{0}}\leq
\theta_{n_{j{1}}}\dots\theta_{n_{j{r}}}\left\|  x\right\|  _{\mathcal{F}_{0}%
}\leq\theta_{1}^{r-1}\theta_{m_{k_{j}-1}}\left\|  x\right\|  _{\mathcal{F}%
_{0}}\leq\theta_{1}^{r-1}\theta_{m_{j}}\left\|  x\right\|  _{\mathcal{F}_{0}%
}.
\]
Therefore
\begin{align*}
\sum_{F\in\mathcal{L}_{r}}t\left(  F\right)  \left\|  Fx\right\|
_{\mathcal{F}_{0}}  &  \leq\sum_{\left\{  j:j<k_{j}\right\}  }\theta_{1}%
^{r-1}\theta_{m_{j}}\left\|  x\right\|  _{\mathcal{F}_{0}}+\sum_{\left\{
j:k_{j}\leq j\right\}  }t\left(  G_{F_{j}}\right)  \left\|  G_{F_{j}%
}x\right\|  _{\mathcal{F}_{0}}\\
&  \leq\theta_{1}^{r-1}\left\|  x\right\|  _{\mathcal{F}_{0}}\sum
_{j=1}^{\infty}\theta_{m_{j}}+\sum_{F\in\mathcal{L}_{r}}t\left(  G_{F}\right)
\left\|  G_{F}x\right\|  _{\mathcal{F}_{0}}.
\end{align*}
Finally,%

\begin{align*}
\left\|  x\right\|   &  =\mathcal{T}^{\mathcal{F}}x=\sum_{r=1}^{\infty}%
\sum_{F\in\mathcal{L}_{r}}t\left(  F\right)  \left\|  Fx\right\|
_{\mathcal{F}_{0}}\\
&  \leq\left\|  x\right\|  _{\mathcal{F}_{0}}\sum_{r=1}^{\infty}\theta
_{1}^{r-1}\sum_{j=1}^{\infty}\theta_{m_{j}}+\sum_{r=1}^{\infty}\sum
_{F\in\mathcal{L}_{r}}t\left(  G_{F}\right)  \left\|  G_{F}x\right\|
_{\mathcal{F}_{0}}\\
&  \leq\frac{2\theta_{m_{1}}}{1-\theta_{1}}\left|  \left|  \left|  x\right|
\right|  \right|  +\left|  \left|  \left|  x\right|  \right|  \right|
=\left(  \frac{2\theta_{m_{1}}}{1-\theta_{1}}+1\right)  \left|  \left|
\left|  x\right|  \right|  \right|  .
\end{align*}
\bigskip
\end{proof}

If $\mathcal{F}$ is a closed subset of $[\mathbb{N}]^{<\infty}$, let
$\mathcal{F}^{\prime}$ be the set of all limit points of $\mathcal{F}$. Define
a transfinite sequence of sets $(\mathcal{F}^{(\alpha)})_{\alpha<\omega_{1}}$
as follows: $\mathcal{F}^{(0)} = \mathcal{F}$, $\mathcal{F}^{(\alpha+1)} =
(\mathcal{F}^{(\alpha)})^{\prime}$ for all $\alpha< \omega_{1}$;
$\mathcal{F}^{(\alpha)} = \cap_{\beta<\alpha}\mathcal{F}^{(\beta)}$ if
$\alpha$ is a countable limit ordinal. If $\mathcal{F}$ is regular, we let
$\iota(\mathcal{F})$ be the unique ordinal $\alpha$ such that $\mathcal{F}%
^{(\alpha)} = \{\emptyset\}$. It is well known that $\iota(\mathcal{S}%
_{\gamma}) = \omega^{\gamma}$ for all $\gamma< \omega_{1}$ \cite[Proposition
4.10]{AA}. The same is true if $\mathcal{S}_{\gamma}$ is replaced by any
$\mathcal{S}^{g}_{\gamma}$.

From now on, we fix a sequence of regular families $(\mathcal{F}_{n})^{\infty
}_{n=0}$ such that $\mathcal{S}_{0} \subseteq\mathcal{F}_{0}$, and a
nonincreasing null sequence $(\theta_{n})^{\infty}_{n=1}$ in $(0,1)$. Denote
the mixed Tsirelson space $T\left(  \mathcal{F}_{0},\left(  \theta
_{n},\mathcal{F}_{n}\right)  _{n=1}^{\infty}\right)  $ by $X$. Let $\alpha_{n}
= \iota(\mathcal{F}_{n})$, $n \in\mathbb{N}\cup\{0\}$. There is no loss of
generality in assuming that $\mathcal{\alpha}_{n}>1$ for all $n\in\mathbb{N}.$
Since $T\left(  \mathcal{F}_{0},\left(  \theta_{n},\mathcal{F}_{n}\right)
_{n=1}^{\infty}\right)  $ is obviously isometric to $T\left(  \mathcal{F}%
_{0},\left(  \theta_{n},\cup_{k=1}^{n}\mathcal{F}_{k}\right)  _{n=1}^{\infty
}\right)  $ via the formal identity, we may also assume that $\left(
\alpha_{n}\right)  _{n=1}^{\infty}$ is a nondecreasing sequence. In the
notation of Proposition \ref{P0}, $\iota\left(  \mathcal{G}_{n}\right)
=\iota\left(  \mathcal{F}_{n}\right)  =\alpha_{n},$ $n\in\mathbb{N}.$ It is
straightforward to check that $\cup_{n=0}^{\infty}\,\mathcal{G}_{n}$ is a
regular family. Relabelling each $\mathcal{G}_{n}$ as $\mathcal{F}_{n},$
$n\in\mathbb{N},$ we may henceforth assume that $\mathcal{S}_{0}%
\subseteq\mathcal{F}_{n}$ for all $n\in\mathbb{N}$ and that $\mathcal{F}%
=\cup_{n=0}^{\infty}\mathcal{F}_{n}$ is regular. Denote $\sup\limits_{n\in
\mathbb{N}}\alpha_{n}$ by $\alpha.$ Note that $\iota\left(  \cup_{n=0}%
^{\infty}\mathcal{F}_{n}\right)  =\alpha\vee\alpha_{0}.$

\section{An estimate on the norm}

\begin{lemma}
\label{L2}Let $\mathcal{G}$ and $\mathcal{H}$ be regular families. Suppose
$\bigcup_{j=1}^{k}F_{j}\in\mathcal{G}\left[  \mathcal{H}\right]  ,$ where
$F_{1}<F_{2}<\dots<F_{k}$. If $F_{j}\notin\mathcal{H}$ for all $j$, $1\leq
j\leq k,$ then $\left\{  \min F_{1},\dots,\min F_{k}\right\}  \in\mathcal{G}.$
\end{lemma}

\begin{proof}
For any nonempty set $G\in\mathcal{G}\left[  \mathcal{H}\right]  ,$ let
$\mathcal{H}\left(  G\right)  =G\cap\left[  1,n\right]  ,$ where $n$ is the
largest integer in $G$ such that $G\cap\left[  1,n\right]  \in\mathcal{H}.$
There is a unique decomposition $G=\bigcup_{j=1}^{k}G_{j},$ where $G_{1}%
,\dots,G_{k}\neq\emptyset$ and $G_{1}=\mathcal{H}\left(  G\right)  ,$
$G_{j+1}=\mathcal{H}\left(  G\setminus\left(  G_{1}\cup\dots\cup G_{j}\right)
\right)  ,$ $1\leq j< k.$ We claim that $\left\{  \min G_{1},\dots,\min
G_{k}\right\}  \in\mathcal{G}.$ To see this, note that since $G\in
\mathcal{G}\left[  \mathcal{H}\right]  ,$ we can write $G=\bigcup_{i=1}^{\ell
}H_{i},$ where $H_{1}<\dots<H_{\ell},$ $H_{1},\dots,H_{\ell}\in\mathcal{H},$
and $\left\{  \min H_{1},\dots,\min H_{\ell}\right\}  \in\mathcal{G}.$
Clearly, $H_{1}\subseteq G_{1}.$ If $k\geq2,$ then $\min H_{2}\leq\min G_{2}.$
If $\max H_{2}>\max G_{2},$ then $G_{2}\subsetneqq H_{2}\subseteq G.$ In
particular, $G_{3}\neq\emptyset$ and $\min G_{3}\in H_{2}.$ Therefore,
$G_{2}\cup\left\{  \min G_{3}\right\}  \in\mathcal{H},$ contrary to the fact
that $G_{2}=\mathcal{H}\left(  G\setminus G_{1}\right)  .$ Thus \thinspace
$\max H_{2}\leq\max G_{2}.$ Continuing this argument, we conclude that $\max
H_{r}\leq\max G_{r}$ for all $1\leq r\leq k.$ It follows that $\left\{  \min
G_{1},\dots,\min G_{k}\right\}  $ is a spreading of $\left\{  \min H_{1}%
,\dots,\min H_{k}\right\}  \in\mathcal{G}.$ Hence $\left\{  \min G_{1}%
,\dots,\min G_{k}\right\}  \in\mathcal{G}.$

Now suppose that $F_{1},\dots,F_{k}$ are as in the statement of the lemma. Let
$G_{j}=\mathcal{H}\left(  F_{j}\right)  ,$ $1\leq j\leq k,$ and let
$G=G_{1}\cup\dots\cup G_{k}.$ Since $G_{j}\subseteq F_{j}$ for $1\leq j\leq
k,$ $G\subseteq\bigcup_{j=1}^{k}F_{j}\in\mathcal{G}\left[  \mathcal{H}\right]
.$ Note that $F_{j}\notin\mathcal{H}$ implies $G_{j}\subsetneqq F_{j}$.
Therefore, $\mathcal{H}\left(  G\right)  =G_{1}$ and
\[
\mathcal{H}\left(  G\setminus\left(  G_{1}\cup\dots\cup G_{j}\right)  \right)
=G_{j+1},\text{ }1\leq j\leq k.
\]
From the previous paragraph, we conclude that $\left\{  \min G_{1},\dots,\min
G_{k}\right\}  \in\mathcal{G}.$ Hence $\left\{  \min F_{1},\dots,\min
F_{k}\right\}  =\left\{  \min G_{1},\dots,\min G_{k}\right\}  \in\mathcal{G}.$
\end{proof}

\begin{proposition}
\label{comparison} \label{P1}Suppose $\varepsilon>0$ and $\mathcal{G}$ is a
regular family. Assume that there exists $m_{0}\in\mathbb{N}$ such that for
all $m\geq m_{0},$ there exist $n_{1},\dots,n_{s}\in\mathbb{N}$ such that
$\theta_{m}<\varepsilon\theta_{n_{1}}\dots\theta_{n_{s}}$ and $\mathcal{F}%
_{m}\subseteq\left[  \mathcal{G},\mathcal{F}_{n_{1}},\dots,\mathcal{F}_{n_{s}%
}\right]  .$ Then there exists a constant $K=K\left(  \varepsilon
,m_{0}\right)  <\infty$ such that for any normalized block basic sequence
$\left(  x_{k}\right)  _{k=1}^{p}$ in $X$ and any real sequence $\left(
a_{k}\right)  _{k=1}^{p},$%
\begin{align}
\left\|  \sum_{k=1}^{p}a_{k}x_{k}\right\|   &  \leq K\left\|  \sum_{k=1}%
^{p}a_{k}e_{i_{k}}\right\|  +2\varepsilon\rho_{1}\left(  \sum_{k=1}^{p}%
a_{k}e_{i_{k}}\right) \label{I1}\\
&  +2\rho_{2}\left(  \sum_{k=1}^{p}a_{k}e_{i_{k}}\right)  +2\varepsilon
\sum_{k=1}^{p}\left|  a_{k}\right|  ,\nonumber
\end{align}
where $i_{k}=\max$\emph{\thinspace supp} $x_{k},$ $1\leq k\leq p$, and
$\rho_{1}$ and $\rho_{2}$ are the norms on the mixed Tsirelson spaces
$T\left(  \mathcal{F},\left(  \theta_{n},\mathcal{F}_{n}\right)
_{n=1}^{\infty}\right)  $ and $T\left(  \mathcal{G},\left(  \theta
_{n},\mathcal{F}_{n}\right)  _{n=1}^{\infty}\right)  $ respectively$.$
\end{proposition}

\begin{proof}
With the given notation, let $x=\sum_{k=1}^{p}a_{k}x_{k}$ and $y=\sum
_{k=1}^{p}a_{k}e_{i_{k}}.$ Also let $G_{k}$ be the integer interval
$(i_{k-1},i_{k}]$ $\left(  i_{0}=0\right)  .$ Since $x\in c_{00},$ there
exists a complete $\left(  \mathcal{F}_{n}\right)  $-admissible tree
$\mathcal{T}$ such that $\left\|  x\right\|  =\mathcal{T}x.$ Each node
$E\in\mathcal{T}$ may be assumed to be contained in the integer interval
$\left[  1,i_{p}\right]  .$ Call a node $E$ $\emph{long}$ if $E\cap G_{k}%
\neq\emptyset$ for at least two values of $k.$ Otherwise, term the node
\emph{short.} Let $N$ be the smallest number such that $\theta_{N}%
\leq\varepsilon.$ Take $\mathcal{E}_{1}$ to be the collection of all minimal
elements in the set of all long nodes $E\in\mathcal{T}$ such that $n\left(
E\right)  > N.$ Minimality is taken with respect to the order (reverse
inclusion) in the tree $\mathcal{T}.$ Similarly, let $\mathcal{E}_{2}$ be the
collection of all minimal elements of the set of all short nodes that are not
in $\cup\left\{  \mathcal{T}_{E}:E\in\mathcal{E}_{1}\right\}  .$ Then let
$\mathcal{E}_{3}$ be the set of all leaves in $\mathcal{T}$ that are not in
$\cup\left\{  \mathcal{T}_{E}:E\in\mathcal{E}_{1}\cup\mathcal{E}_{2}\right\}
.$ Observe that%
\[
\mathcal{T}x\leq\sum_{j=1}^{3}\sum_{E\in\mathcal{E}_{j}}t\left(  E\right)
\left\|  Ex\right\|
\]
The proof of the proposition is completed by combining Lemmas \ref{L4},
\ref{L5}, \ref{L6}, and \ref{L7} below.
\end{proof}

\begin{lemma}
\label{L4}$\sum\limits_{E\in\mathcal{E}_{1}}t\left(  E\right)  \left\|
Ex\right\|  \leq$ $2\varepsilon\sum\limits_{k=1}^{p}\left|  a_{k}\right|  .$
\end{lemma}

\begin{proof}
Arrange the nodes in $\mathcal{E}_{1}$ from left to right as $E_{1}%
<\dots<E_{r}.$ Since $n\left(  E_{j}\right)  >N,$ $t\left(  E_{j}\right)
<\theta_{N}\leq\varepsilon.$ For $1\leq j\leq r,$ let $J_{j}=\left\{
k:G_{k}\cap E_{j}\neq\emptyset\right\}  .$ Then $J_{1}\leq\dots\leq J_{r},$
and $\left|  J_{j}\right|  \geq2$ for all $j.$ Hence $\sum_{j=1}^{r}\sum_{k\in
J_{j}}\left|  a_{k}\right|  \leq2\sum_{k=1}^{p}\left|  a_{k}\right|  .$ It
follows that
\[
\sum_{E\in\mathcal{E}_{1}}t\left(  E\right)  \left\|  Ex\right\|  \leq
\sum_{j=1}^{r}t\left(  E_{j}\right)  \sum_{k\in J_{j}}\left|  a_{k}\right|
\leq\varepsilon\sum_{j=1}^{r}\sum_{k\in J_{j}}\left|  a_{k}\right|
\leq2\varepsilon\sum_{k=1}^{p}\left|  a_{k}\right|  .
\]
\end{proof}

\begin{lemma}
\label{L5}$\sum_{E\in\mathcal{E}_{3}}t\left(  E\right)  \left\|  Ex\right\|
\leq2\left\|  \sum_{k=1}^{p}a_{k}e_{i_{k}}\right\|  .$
\end{lemma}

\begin{proof}
Since any node $E\in\mathcal{E}_{3}$ is a leaf in the complete tree
$\mathcal{T}$ for $x,$ $\left\|  Ex\right\|  =\left\|  Ex\right\|
_{\mathcal{F}_{0}}.$ Choose $E_{0}\in\mathcal{F}_{0}$ such that $E_{0}%
\subseteq E$ and $\left\|  Ex\right\|  = \|E_{0}x\| = \left\|  E_{0}x\right\|
_{\ell^{1}}.$ Let $J_{E}=\left\{  k:G_{k}\cap E_{0}\neq\emptyset\right\}  .$
For each $k\in J_{E},$ choose $j_{k}\in G_{k}\cap E_{0}$ and set $z=\sum
_{E\in\mathcal{E}_{3}}\sum_{k\in J_{E}}a_{k}e_{j_{k}}.$ Because each
$E\in\mathcal{E}_{3}$ is a long node, each $k$ belongs to at most two $J_{E}.$
It follows that $\left\|  z\right\|  \leq2\left\|  \sum_{k=1}^{p}a_{k}%
e_{i_{k}}\right\|  .$ Now%
\begin{align*}
\sum_{E\in\mathcal{E}_{3}}t\left(  E\right)  \left\|  Ex\right\|   &  =
\sum_{E\in\mathcal{E}_{3}}t\left(  E\right)  \left\|  E_{0}x\right\|  \leq
\sum_{E\in\mathcal{E}_{3}}t\left(  E\right)  \sum_{k\in J_{E}}\left|
a_{k}\right| \\
&  \leq\sum_{E\in\mathcal{E}_{3}}t\left(  E\right)  \left\|  E_{0}z\right\|
_{\ell^{1}} \leq\sum_{E\in\mathcal{E}_{3}}t\left(  E\right)  \left\|
Ez\right\|  _{\mathcal{F}_{0}}\\
&  \leq\left\|  z\right\|  \leq2\bigl\| \sum_{k=1}^{p}a_{k}e_{i_{k}}\bigr\| .
\end{align*}
\end{proof}

Observe that any ancestor $F$ of any node in $\mathcal{E}_{2}$ must be a long
node such that $n(F)\leq N$. Subdivide $\mathcal{E}_{2}$ into two parts
$\mathcal{E}_{21}$ and $\mathcal{E}_{22}$ according to whether the node $E$ in
question satisfies $n\left(  E\right)  >N$ or $n\left(  E\right)  \leq N$.

\begin{lemma}
\label{L6}$\sum\limits_{E\in\mathcal{E}_{21}}t\left(  E\right)  \left\|
Ex\right\|  \leq\frac{2}{\theta_{m_{0}}}\left\|  y\right\|  +2\varepsilon
\rho_{1}\left(  y\right)  +2\rho_{2}\left(  y\right)  .$
\end{lemma}

\begin{proof}
Let $\mathcal{D}$ be the set of all nodes that are immediate predecessors of
some node in $\mathcal{E}_{21}$. Let us first show that any two distinct nodes
$D$ and $D^{\prime}$ in $\mathcal{D}$ are mutually incomparable. Indeed,
suppose that $D$ is an ancestor of $D^{\prime}$. Let $E$ and $E^{\prime}$ be
immediate successors of $D$ and $D^{\prime}$ respectively that are in
$\mathcal{E}_{21}$. Consider the immediate successor $D^{\prime\prime}$ of $D$
such that $D^{\prime\prime}\supseteq D^{\prime}$. Since $D^{\prime\prime}$ and
$E$ are both immediate successors of $D$, $n\left(  D^{\prime\prime}\right)
=n\left(  E\right)  $. But $n\left(  D^{\prime\prime}\right)  \leq N$ since
$D^{\prime\prime}$ is an ancestor of $E^{\prime}\in\mathcal{E}_{21}$, while
$n\left(  E\right)  >N$ by definition of $\mathcal{E}_{21}.$ Thus $D$ and
$D^{\prime}$ must be mutually incomparable. List the elements in $\mathcal{D}$
from left to right as $D_{1}<D_{2}<\dots<D_{r}$. If $1\leq j\leq r$ and $1\leq
k\leq p,$ let $\mathcal{D}_{jk}=\left\{  E\in\mathcal{E}_{21}:E\subseteq
D_{j}\cap G_{k}\right\}  $ and $J_{j}=\left\{  k:\mathcal{D}_{jk}\neq
\emptyset\right\}  .$ By the preceding argument, each $E$ in $\bigcup_{k\in
J_{j}}\mathcal{D}_{jk}$ is an immediate successor of $D_{j}$. Given $k\in
J_{j},$ choose $E_{jk}\in\mathcal{D}_{jk}$ and $\ell_{jk}\in E_{jk}.$ As in
the proof of Lemma \ref{L5}, note that each $k$ belongs to at most two $J_{j}$
because each $D_{j}$ is a long node$.$ Hence $\left\|  w\right\|
\leq2\left\|  y\right\|  $ and $\rho_{i}\left(  w\right)  \leq2\rho_{i}\left(
y\right)  ,$ $i=1,2,$ where $w=\sum_{j=1}^{r}\sum_{k\in J_{j}}a_{k}%
e_{\ell_{jk}}.$ For each $j,$ let $m=m\left(  j\right)  $ be the common value
of $n(E)$ for all $E \in\bigcup_{k\in J_{j}}\mathcal{D}_{jk}$. In particular,
$\bigcup_{k\in J_{j}}\mathcal{D}_{jk}$ is $\mathcal{F}_{m}$-admissible.
Consider the set $M=\left\{  j:m\left(  j\right)  <m_{0}\right\}  .$ If $j\in
M,$ then%
\begin{align*}
\sum_{k\in J_{j}}\sum_{E\in\mathcal{D}_{jk}}t\left(  E\right)  \left\|
Ex\right\|   &  =\sum_{k\in J_{j}}t\left(  D_{j}\right)  \theta_{m}\sum
_{E\in\mathcal{D}_{jk}}\left\|  Ex\right\|  \leq t\left(  D_{j}\right)
\sum_{k\in J_{j}}\left|  a_{k}\right| \\
&  \leq\frac{t\left(  D_{j}\right)  }{\theta_{m_{0}}}\theta_{m}\sum_{k\in
J_{j}}\sum_{E\in\mathcal{D}_{jk}}\left\|  E\left(  D_{j}w\right)  \right\|
_{\mathcal{S}_{0}}\\
&  \leq\frac{t\left(  D_{j}\right)  }{\theta_{m_{0}}}\left\|  D_{j}w\right\|
.
\end{align*}
Hence
\begin{align}
\label{E0}\sum_{j\in M}\sum_{k\in J_{j}}\sum_{E\in\mathcal{D}_{jk}}t\left(
E\right)  \left\|  Ex\right\|   &  \leq\frac{1}{\theta_{m_{0}}}\sum_{j\in
M}t\left(  D_{j}\right)  \left\|  D_{j}w\right\| \\
&  \leq\frac{1}{\theta_{m_{0}}}\left\|  w\right\|  \leq\frac{2}{\theta_{m_{0}%
}}\left\|  y\right\|  .\nonumber
\end{align}
If $j\notin M,$ choose $n_{1},\dots,n_{s}\in\mathbb{N}$ as in the hypothesis
of Proposition \ref{P1}. Note that $I_{j}=\left\{  \ell_{jk}:k\in
J_{j}\right\}  \in\mathcal{F}_{m}.$ Partition $J_{j}$ into $J_{j}^{\prime}$
and $J_{j}^{\prime\prime}$ so that $J_{j}^{\prime}$ consists of all $k\in
J_{j}$ such that $\mathcal{D}_{jk}\text{ is }\left[  \mathcal{F}_{n_{1}}%
,\dots,\mathcal{F}_{n_{s}}\right]  \text{-admissible} $ and $J_{j}%
^{\prime\prime}=J_{j}\setminus J_{j}^{\prime}.$ Set $I_{j}^{\prime}=\{
\ell_{jk}:k\in J_{j}^{\prime}\} .$ Then%
\begin{align*}
\sum_{k\in J_{j}^{\prime}}\sum_{E\in\mathcal{D}_{jk}}t\left(  E\right)
\left\|  Ex\right\|   &  =t\left(  D_{j}\right)  \sum_{k\in J_{j}^{\prime}%
}\theta_{m}\sum_{E\in\mathcal{D}_{jk}}\left\|  Ex\right\| \\
&  \leq\varepsilon t\left(  D_{j}\right)  \sum_{k\in J_{j}^{\prime}}%
\theta_{n_{1}}\dots\theta_{n_{s}}\sum_{E\in\mathcal{D}_{jk}}\left\|
Ex\right\| \\
&  \leq\varepsilon t\left(  D_{j}\right)  \sum_{k\in J_{j}^{\prime}}\left|
a_{k}\right|  \leq\varepsilon t\left(  D_{j}\right)  \left\|  I_{j}^{\prime
}(D_{j}w)\right\|  _{\ell^{1}}\\
&  \leq\varepsilon t\left(  D_{j}\right)  \left\|  D_{j}w\right\|
_{\mathcal{F}_{m}}\leq\varepsilon t\left(  D_{j}\right)  \left\|
D_{j}w\right\|  _{\mathcal{F}}.
\end{align*}
Hence
\begin{align}
\label{E1}\sum_{j\notin M}\sum_{k\in J_{j}^{\prime}}\sum_{E\in\mathcal{D}%
_{jk}}t\left(  E\right)  \left\|  Ex\right\|   &  \leq\varepsilon\sum
_{j=1}^{r}t\left(  D_{j}\right)  \left\|  D_{j}w\right\|  _{\mathcal{F}}\\
&  \leq\varepsilon\rho_{1}\left(  w\right)  \leq2\varepsilon\rho_{1}\left(
y\right)  .\nonumber
\end{align}
On the other hand, since $\bigcup_{k\in J_{j}^{\prime\prime}}\mathcal{D}_{jk}$
is $\mathcal{F}_{m}$- and thus $\left[  \mathcal{G},\mathcal{F}_{n_{1}}%
,\dots,\mathcal{F}_{n_{s}}\right]  $-ad\-miss\-ible, while $\mathcal{D}_{jk}$
is not $\left[  \mathcal{F}_{n_{1}},\dots,\mathcal{F}_{n_{s}}\right]
$-admissible for all $k\in J_{j}^{\prime\prime},$
\[
\bigl\{  \min\cup_{E\in\mathcal{D}_{jk}}E:k\in J_{j}^{\prime\prime}\bigr
\}  \in\mathcal{G}%
\]
by Lemma \ref{L2}. Thus $I_{j}^{\prime\prime}=\bigl\{  \ell_{jk}:k\in
J_{j}^{\prime\prime}\bigr\}  \in\mathcal{G}.$ Consequently,
\begin{align*}
\sum_{k\in J_{j}^{\prime\prime}}\sum_{E\in\mathcal{D}_{jk}}t\left(  E\right)
\left\|  Ex\right\|   &  =t\left(  D_{j}\right)  \sum_{k\in J_{j}%
^{\prime\prime}}\theta_{m}\sum_{E\in\mathcal{D}_{jk}}\left\|  Ex\right\|  \leq
t\left(  D_{j}\right)  \sum_{k\in J_{j}^{\prime\prime}}\left|  a_{k}\right| \\
&  \leq t\left(  D_{j}\right)  \left\|  I_{j}^{\prime\prime}(D_{j}w)\right\|
_{\ell^{1}} \leq t\left(  D_{j}\right)  \left\|  D_{j}w\right\|
_{\mathcal{G}}.
\end{align*}
Therefore
\begin{equation}
\sum_{j\notin M}\sum_{k\in J_{j}^{\prime\prime}}\sum_{E\in\mathcal{D}_{jk}%
}t\left(  E\right)  \left\|  Ex\right\|  \leq\sum_{j=1}^{r}t\left(
D_{j}\right)  \left\|  D_{j}w\right\|  _{\mathcal{G}}\leq\rho_{2}\left(
w\right)  \leq2\rho_{2}\left(  y\right)  . \label{E2}%
\end{equation}
Combining inequalities (\ref{E0}), (\ref{E1}) and (\ref{E2}) completes the proof.
\end{proof}

\begin{lemma}
\label{L7}$\sum\limits_{E\in\mathcal{E}_{22}}t\left(  E\right)  \left\|
Ex\right\|  \leq\frac{2}{\theta_{N}}\left\|  y\right\|  $.
\end{lemma}

\begin{proof}
For $1\leq k\leq p,$ let $\mathcal{E}_{22}\left(  k\right)  =\left\{
E\in\mathcal{E}_{22}:E\subseteq G_{k}\right\}  .$ If $\mathcal{E}_{22}\left(
k\right)  \neq\emptyset,$ denote by $\mathcal{P}_{k}$ the collection of all
minimal elements in the set of all nodes that are immediate predecessors of
some node in $\mathcal{E}_{22}\left(  k\right)  .$ Observe that if
$P\in\mathcal{P}_{k},$ then $P$ is a long node and $P\cap G_{k}\neq\emptyset.$
Hence $\left|  \mathcal{P}_{k}\right|  \leq2.$ For each $P\in\mathcal{P}_{k},$
choose an immediate successor $E_{P}$ of $P$ such that $E_{P}\in
\mathcal{E}_{22}\left(  k\right)  ,$ then fix $j_{P}\in E_{P}.$ Note that the
nodes in $\{E_{P}: P \in\cup^{p}_{k=1}\mathcal{P}_{k}\}$ are pairwise
disjoint. Set $v=\sum_{k=1}^{p}a_{k}\sum_{P\in\mathcal{P}_{k}}e_{j_{P}}.$
Since $\left|  \mathcal{P}_{k}\right|  \leq2,$ $\left\|  v\right\|
\leq2\left\|  y\right\|  .$ Notice that $t\left(  E_{P}\right)  =\theta
_{n\left(  E_{P}\right)  }t\left(  P\right)  \geq\theta_{N}t\left(  P\right)
$ since $E\in\mathcal{E}_{22}$ implies $n\left(  E_{P}\right)  \leq N.$ Now%
\begin{align*}
\sum_{E\in\mathcal{E}_{22}}t\left(  E\right)  \left\|  Ex\right\|   &
=\sum_{k=1}^{p}\sum_{E\in\mathcal{E}_{22}\left(  k\right)  }t\left(  E\right)
\left\|  Ex\right\|  =\sum_{k=1}^{p}\sum_{P\in\mathcal{P}_{k}}\underset
{E\subseteq P}{\sum_{E\in\mathcal{E}_{22}\left(  k\right)  }}t\left(
E\right)  \left\|  Ex\right\| \\
&  \leq\sum_{k=1}^{p}\sum_{P\in\mathcal{P}_{k}}t\left(  P\right)  \left\|
P(G_{k}x)\right\|  \leq\sum_{k=1}^{p}\sum_{P\in\mathcal{P}_{k}}t\left(
P\right)  \left|  a_{k}\right| \\
&  =\sum_{k=1}^{p}\sum_{P\in\mathcal{P}_{k}}t\left(  P\right)  \left\|
E_{P}v\right\|  _{\mathcal{S}_{0}}\leq\frac{1}{\theta_{N}}\sum_{k=1}^{p}%
\sum_{P\in\mathcal{P}_{k}}t\left(  E_{P}\right)  \left\|  E_{P}v\right\|
_{\mathcal{F}_{0}}\\
&  \leq\frac{1}{\theta_{N}}\left\|  v\right\|  \leq\frac{2}{\theta_{N}%
}\left\|  y\right\|  .
\end{align*}
\end{proof}

Observe that in the preceding proof, the hypothesis of Proposition \ref{P1}
(that is, the existence of the family $\mathcal{G}$) is used only in Lemma 6.
One may consider mixed Tsirelson spaces $Z=T( \mathcal{F}_{0},\left(
\theta_{n},\mathcal{F}_{n}) _{n=1}^{\ell}\right)  $ determined by finitely
many regular families, defined in the obvious way. For such spaces, it is
worthwhile to observe the following corollary of the proof of Proposition
\ref{P1}.

\begin{corollary}
\label{equivalence} Let the space $Z$ be as above. There exists a constant
$K<\infty$ such that for any normalized block basic sequence $\left(
x_{k}\right)  _{k=1}^{p}$ in $Z$ and any $\left(  a_{k}\right)  \in c_{00},$%
\[
\frac{1}{2}\left\|  \sum_{k=1}^{p}a_{k}e_{i_{k}}\right\|  \leq\left\|
\sum_{k=1}^{p}a_{k}x_{k}\right\|  \leq K\left\|  \sum_{k=1}^{p}a_{k}e_{i_{k}%
}\right\|  ,
\]
where $i_{k}=\max$\emph{\thinspace supp} $x_{k},$ $1\leq k\leq p.$
\end{corollary}

\begin{proof}
If $j_{k}=\min\,$\thinspace supp $x_{k},$ $1\leq k\leq p,$ then
\[
\frac{1}{2}\left\|  \sum_{k=1}^{p}a_{k}e_{i_{k}}\right\|  \leq\left\|
\sum_{k=1}^{p}a_{k}e_{j_{k}}\right\|  \leq\left\|  \sum_{k=1}^{p}a_{k}%
x_{k}\right\|  .
\]
On the other hand, in the notation of the proof of Proposition \ref{P1}, take
$N=\ell$.
Then $\mathcal{E}_{1}= \mathcal{E}_{21} = \emptyset$. In particular, the
hypothesis in Proposition \ref{comparison} is no longer required since Lemma
\ref{L6} is not needed any more. Lemmas \ref{L5} and \ref{L7} give the desired result.
\end{proof}

\section{Bounds on the $\ell^{1}$-index}

Let us recall the relevant terminology concerning trees. A \emph{tree} on a
set $S$ is a subset $T$ of $\cup_{n=1}^{\infty}S^{n}$ such that $(x_{1}%
,\dots,x_{n})\in T$ whenever $n\in\mathbb{N}$ and $(x_{1},\dots,x_{n+1})\in
T$. If $(x_{1},\dots,x_{n}) \in T$ and $1 \leq m < n$, the sequence
$(x_{1},\dots, x_{m})$ is said to be an ancestor of $(x_{1},\dots,x_{n})$. A
tree $T$ is \emph{well-founded} if there is no infinite sequence $(x_{n})$ in
$S$ such that $(x_{1},\dots,x_{n})\in T$ for all $n$. Given a well-founded
tree $T$, we define the \emph{derived tree} $D(T)$ to be the set of all
$(x_{1},\dots,x_{n})\in T$ such that $(x_{1},\dots,x_{n},x)\in T$ for some
$x\in S$. Inductively, we let $D^{0}(T)=T$, $D^{\alpha+1}(T)=D(D^{\alpha}%
(T))$, and $D^{\alpha}(T)=\cap_{\beta<\alpha}D^{\beta}(T)$ if $\alpha$ is a
limit ordinal. The \emph{order} of a well-founded tree $T$ is the smallest
ordinal $o(T)$ such that $D^{o(T)}(T)=\emptyset$. If $E$ is a Banach space and
$1\leq K<\infty$, an $\ell^{1}$-$K$ tree on $E$ is a tree $T$ on $S(E)=\{x\in
E:\Vert x\Vert=1\}$ such that $\Vert\sum_{i=1}^{n}a_{i}x_{i}\Vert\geq
K^{-1}\sum_{i=1}^{n}|a_{i}|$ whenever $(x_{1},\dots,x_{n})\in T$ and
$(a_{i})\subseteq\mathbb{R}$. If $E$ has a basis $(e_{i})$, a \emph{block
tree} on $E$ is a tree $T$ on $E$ so that every $(x_{1},\dots,x_{n})\in T$ is
a finite block basis of $(e_{i})$. An $\ell^{1}$-$K$-block tree on $E$ is a
block tree that is also an $\ell^{1}$-$K$ tree. The index $I(E,K)$ is defined
to be $\sup\{o(T):T\text{ is an $\ell^{1}$-$K$ tree on $E$}\}$. If $E$ has a
basis $(e_{i})$, the index $I_{b}(E,K)$ is defined similarly, with the
supremum taken over all $\ell^{1}$-$K$ block trees. The Bourgain $\ell^{1}%
$-\emph{index} of $E$ is the ordinal $I(E)=\sup\{I(E,K):1\leq K<\infty\}$. The
index $I_{b}(E)$ is defined similarly. Bourgain proved that if $E$ is a
separable Banach space not containing a copy of $\ell^{1}$, then
$I(E)<\omega_{1}$ \cite{B}. Judd and Odell \cite{JO} showed that $I(E)$ and
$I_{b}(E)$ are closely related for a Banach space $E$ with a basis. Precisely,
if $I_{b}(E) = \omega^{n}$ for some $n < \omega$, then $I(E) = \omega^{n}$ or
$\omega^{n+1}$, while $I_{b}(E) = I(E)$ if $I_{b}(E) \geq\omega^{\omega}$. We
refer the reader to \cite{JO} and \cite{AJO} for in depth discussions of these
and related indices.

Our concern for the rest of the paper is the calculation of the index
$I_{b}(X)$, where $X$ is the mixed Tsirelson space $T\left(  \mathcal{F}%
_{0},\left(  \theta_{n},\mathcal{F}_{n}\right)  _{n=1}^{\infty}\right)  $. We
begin with an easy lower bound on $I_{b}\left(  X\right)  .$

\begin{proposition}
\label{P8}$I_{b}\left(  X\right)  \geq\alpha_{0}\cdot\sup\limits_{n\in
\mathbb{N}}\alpha_{n}^{\omega}.$
\end{proposition}

\begin{proof}
For all $m,n\in\mathbb{N},$ denote the family $\overset{m\text{ times}%
}{\overbrace{[\mathcal{F}_{n},\dots,\mathcal{F}_{n}}},\mathcal{F}_{0}]$ by
$\mathcal{B}_{mn}.$ Observe that $\iota\left(  \mathcal{B}_{mn}\right)
=\alpha_{0}\cdot\alpha_{n}^{m}$ for all $m,n\in\mathbb{N}$ by
\cite[Proposition 10]{LT}. For any $\left(  a_{k}\right)  \in c_{00},$
$\left\|  \sum a_{k}e_{k}\right\|  \geq\theta_{n}^{m}\left\|  \sum a_{k}%
e_{k}\right\|  _{\mathcal{B}_{mn}}.$ Thus $I_{b}\left(  X,\theta_{n}%
^{m}\right)  \geq\iota\left(  \mathcal{B}_{mn}\right)  =\alpha_{0}\cdot
\alpha_{n}^{m}$ for all $m,n\in\mathbb{N}.$ Therefore,%
\[
I_{b}\left(  X\right)  \geq\sup_{m,n\in\mathbb{N}}\alpha_{0}\cdot\alpha
_{n}^{m}=\alpha_{0}\cdot\sup\limits_{n\in\mathbb{N}}\alpha_{n}^{\omega}.
\]
\end{proof}

In the remainder of this section, we apply Proposition \ref{P1} to obtain an
upper bound on the $\ell^{1}$-index of $X.$ For each $n\in\mathbb{N},$ let%
\[
\mathcal{C}\left(  n\right)  =\left\{  \left(  0,n_{1},\dots,n_{s}\right)
:n_{1}, \dots, n_{s}, s\in\mathbb{N},\text{ }n_{1}+n_{2}+\dots+n_{s}\leq
n\right\}
\]
and
\[
\pi_{n}=\sup\left\{  \theta_{n_{1}}\dots\theta_{n_{s}}:n_{1}+\dots
+n_{s}>n\right\}  .
\]
Obviously $\mathcal{C}\left(  n\right)  $ is a finite set. Denote its
cardinality by $p\left(  n\right)  .$ It is clear that $\lim
\limits_{n\rightarrow\infty}\pi_{n}=0.$

\begin{lemma}
\label{P3}Suppose that $\mathcal{H}$ is a regular family containing
$\mathcal{S}_{0}$ and that $\rho$ is the norm on the space $T\left(
\mathcal{H},\left(  \theta_{n},\mathcal{F}_{n}\right)  _{n=1}^{\infty}\right)
.$ For all $x\in c_{00}$ and all $n\in\mathbb{N},$ we have
\[
\rho\left(  x\right)  \leq\pi_{n}\left\|  x\right\|  _{\ell^{1}}+p\left(
n\right)  \left\|  x\right\|  _{\mathcal{M}\left[  \mathcal{H}\right]  },
\]
where $\mathcal{M}=\cup_{\left(  0,n_{1},\dots,n_{s}\right)  \in
\mathcal{C}\left(  n\right)  }\left[  \mathcal{F}_{n_{1}},\dots,\mathcal{F}%
_{n_{s}}\right]  .$
\end{lemma}

\begin{proof}
There exists an $\left(  \mathcal{F}_{n}\right)  $-admissible tree
$\mathcal{T}$ such that
\[
\rho\left(  x\right)  =\sum_{E\in\mathcal{L}}t\left(  E\right)  \left|
\left|  Ex\right|  \right|  _{\mathcal{H}},
\]
where $\mathcal{L}$ is the set of all leaves of $\mathcal{T}.$ Let
$\mathcal{L}_{\left(  n_{1},\dots,n_{s}\right)  }$ be the set of all
$E\in\mathcal{L}$ such that $h\left(  E\right)  =\left(  0,n_{1},\dots
,n_{s}\right)  $. Then%
\[
\rho\left(  x\right)  =\left(  \sum_{\left(  0,n_{1},\dots,n_{s}\right)
\in\mathcal{C}\left(  n\right)  }\,+\sum_{\left(  0,n_{1},\dots,n_{s}\right)
\notin\mathcal{C}\left(  n\right)  }\right)  \,\sum_{E\in\mathcal{L}_{\left(
n_{1},\dots,n_{s}\right)  }}t\left(  E\right)  \left\|  Ex\right\|
_{\mathcal{H}}.
\]
If $\left(  0,n_{1},\dots,n_{s}\right)  \in\mathcal{C}\left(  n\right)  ,$
then $\mathcal{L}_{\left(  n_{1},\dots,n_{s}\right)  }$ is $\left[
\mathcal{F}_{n_{1}},\dots,\mathcal{F}_{n_{s}}\right]  $-admissible and thus
$\mathcal{M}$-admissible. Since $t\left(  E\right)  \leq1$ for all $E$,
\[
\sum_{E\in\mathcal{L}_{\left(  n_{1},\dots,n_{s}\right)  }}t\left(  E\right)
\left\|  Ex\right\|  _{\mathcal{H}} \leq\sum_{E\in\mathcal{L}_{\left(
n_{1},\dots,n_{s}\right)  }}\left\|  Ex\right\|  _{\mathcal{H}}\leq\left\|
x\right\|  _{\mathcal{M}\left[  \mathcal{H}\right]  }.
\]
Therefore,
\[
\sum_{\left(  0,n_{1},\dots,n_{s}\right)  \in\mathcal{C}\left(  n\right)
}\sum_{E\in\mathcal{L}_{\left(  n_{1},\dots,n_{s}\right)  }}t\left(  E\right)
\left\|  Ex\right\|  _{\mathcal{H}}\leq p\left(  n\right)  \left\|  x\right\|
_{\mathcal{M}\left[  \mathcal{H}\right]  }.
\]
On the other hand, since $t(E) = \theta_{n_{1}}\dots\theta_{n_{s}} \leq\pi
_{n}$ if $E \notin\mathcal{L}_{(n_{1},\dots,n_{s})}$,
\begin{align*}
&  \quad\sum_{\left(  0,n_{1},\dots,n_{s}\right)  \notin\mathcal{C}\left(
n\right)  }\sum_{E\in\mathcal{L}_{\left(  n_{1},\dots,n_{s}\right)  }}t\left(
E\right)  \left\|  Ex\right\|  _{\mathcal{H}}\\
&  \leq\sum_{\left(  0,n_{1},\dots,n_{s}\right)  \notin\mathcal{C}\left(
n\right)  }\sum_{E\in\mathcal{L}_{\left(  n_{1},\dots,n_{s}\right)  }}\pi
_{n}\left\|  Ex\right\|  _{\mathcal{H}}\leq\pi_{n}\left\|  x\right\|
_{\ell^{1}}.
\end{align*}
\end{proof}

\begin{lemma}
Let $\mathcal{M}$ be as defined in Lemma \ref{P3}, then $\iota\left(
\mathcal{M}\right)  \leq\alpha_{n}^{n}.$
\end{lemma}

\begin{proof}
The lemma follows immediately from the fact that
\[
\iota\left(  \mathcal{H}\left[  \mathcal{N}\right]  \right)  \leq\iota\left(
\mathcal{N}\right)  \cdot\iota\left(  \mathcal{H}\right)
\]
if $\mathcal{H}$ and $\mathcal{N}$ are regular families of finite subsets of
$\mathbb{N}$ (cf. \cite[Proposition 10]{LT}).
\end{proof}

\begin{proposition}
\emph{(\cite[Proposition 12]{LT})}\label{L8} Let $T$ be a well-founded block
tree on some basis $\left(  e_{i}\right)  .$ Define
\[
\mathcal{H}\left(  T\right)  =\left\{  \left\{  \max\text{\emph{supp\thinspace
}}x_{i}:i=1,\dots,n\right\}  :\left(  x_{1},x_{2},\dots,x_{n}\right)  \in
T\right\}
\]
and
\[
\mathcal{G}\left(  T\right)  =\{G:G\text{ is a spreading of a subset of some
$F\in\mathcal{H}(T)$}\}.
\]
If $\mathcal{G}(T)$ is compact, then $\iota\left(  \mathcal{G}\left(
T\right)  \right)  \geq o\left(  T\right)  .$
\end{proposition}

Given a countable ordinal $\eta,$ define the order (or the logarithm)
$\ell\left(  \eta\right)  $ of the ordinal $\eta$ to be $\gamma_{1},$ where
$\eta=\omega^{\gamma_{1}}\cdot k_{1}+\dots+\omega^{\gamma_{p}}\cdot k_{p}$ in
Cantor normal form. Clearly, $\ell\left(  \eta_{1}\cdot\eta_{2}\right)
=\ell\left(  \eta_{1}\right)  +\ell\left(  \eta_{2}\right)  .$ Therefore
$\ell\left(  \eta^{n}\right)  =\ell\left(  \eta\right)  \cdot n$ and
$\ell\left(  \eta^{\omega}\right)  =\ell\left(  \eta\right)  \cdot\omega.$
Obviously, if $\ell\left(  \eta\right)  = \gamma,$ then $\omega^{\gamma}%
\leq\eta<\omega^{\gamma+1}.$ Observe that in the notation of Proposition
\ref{P1}, if we take $\rho$ to be the norm on the space $T\left(
\mathcal{F}_{0}\cup\mathcal{G},\left(  \theta_{n},\mathcal{F}_{n}\right)
_{n=1}^{\infty}\right)  ,$ then $\left\|  \cdot\right\|  \leq\rho$ and
$\rho_{2}\leq\rho.$ Thus inequality (\ref{I1}) implies%
\[
\left\|  \sum_{k=1}^{p}a_{k}x_{k}\right\|  \leq\left(  K+2\right)  \rho\left(
\sum_{k=1}^{p}a_{k}e_{i_{k}}\right)  +4\varepsilon\sum_{k=1}^{p}\left|
a_{k}\right|  .
\]
If $\left(  x_{k}\right)  _{k=1}^{n}$ and $\left(  y_{k}\right)  _{k=1}^{n}$
are sequences in possibly different normed spaces, and $0<K<\infty,$ we write
$\left(  x_{k}\right)  _{k=1}^{n}\overset{K}{\succeq}\left(  y_{k}\right)
_{k=1}^{n}$ to mean $K\left\|  \sum_{k=1}^{n}a_{k}x_{k}\right\|  \geq\left\|
\sum_{k=1}^{n}a_{k}y_{k}\right\|  $ for all $\left(  a_{k}\right)  \in c_{00}$.

\begin{proposition}
\label{P11}Suppose for all $\varepsilon>0,$ there exist a regular family
$\mathcal{G}_{\varepsilon}$ and $m_{0}\in\mathbb{N}$ such that for all $m\geq
m_{0},$ there exist $n_{1},\dots,n_{s}\in\mathbb{N}$ satisfying $\theta
_{m}<\varepsilon\theta_{n_{1}}\dots\theta_{n_{s}}$ and $\mathcal{F}%
_{m}\subseteq\left[  \mathcal{G}_{\varepsilon},\mathcal{F}_{n_{1}}%
,\dots,\mathcal{F}_{n_{s}}\right]  .$ Then%
\[
I_{b}\left(  X\right)  \leq\sup_{\varepsilon>0}\sup\limits_{n\in\mathbb{N}%
}\left[  \left(  \alpha_{0}\vee\iota\left(  \mathcal{G}_{\varepsilon}\right)
\right)  \cdot\alpha_{n}^{\omega}\right]  .
\]
\end{proposition}

\begin{proof}
Suppose otherwise. There exists $H>1$ and an $\ell^{1}$-$H$-block tree $T$ on
$X$ such that
\[
o\left(  T\right)  >\sup\limits_{\varepsilon>0}\sup\limits_{n\in\mathbb{N}%
}\left[  \left(  \alpha_{0}\vee\iota\left(  \mathcal{G}_{\varepsilon}\right)
\right)  \cdot\alpha_{n}^{\omega}\right]  .
\]
Pick $\varepsilon_{0}<\frac{1}{8H}.$ According to Proposition \ref{P1} and the
remark above, there exists a constant $K$ such that%
\[
\left\|  \sum_{k=1}^{n}a_{k}x_{k}\right\|  \leq K\rho\left(  \sum_{k=1}%
^{n}a_{k}e_{i_{k}}\right)  +4\varepsilon_{0}\sum_{k=1}^{n}\left|
a_{k}\right|
\]
for all $\left(  a_{k}\right)  \in c_{00},$ where $\rho$ is the norm on
$T\left(  \mathcal{F}_{0}\cup\mathcal{G}_{\varepsilon_{0}},\left(  \theta
_{n},\mathcal{F}_{n}\right)  _{n=1}^{\infty}\right)  .$ Let $\ell\left(
\alpha_{n}\right)  =\gamma_{n}$ for all $n\in\mathbb{N}$ and $\ell\left(
\alpha_{0}\vee\iota\left(  \mathcal{G}_{\varepsilon_{0}}\right)  \right)
=\gamma$. Then
\begin{align*}
\ell\left(  \left(  \alpha_{0}\vee\iota\left(  \mathcal{G}_{\varepsilon_{0}%
}\right)  \right)  \cdot\sup\limits_{n\in\mathbb{N}}\alpha_{n}^{\omega
}\right)   &  =\ell\left(  \alpha_{0}\vee\iota\left(  \mathcal{G}%
_{\varepsilon_{0}}\right)  \right)  +\ell\left(  \sup\limits_{n\in\mathbb{N}%
}\alpha_{n}^{\omega}\right) \\
&  \geq\ell\left(  \alpha_{0}\vee\iota\left(  \mathcal{G}_{\varepsilon_{0}%
}\right)  \right)  +\ell\left(  \alpha_{n}^{\omega}\right)  =\gamma+\gamma
_{n}\cdot\omega
\end{align*}
for all $n\in\mathbb{N}.$ Hence $o\left(  T\right)  >\omega^{\gamma+\gamma
_{n}\cdot\omega}$ for all $n\in\mathbb{N}.$ Given $F\in\mathcal{H}\left(
T\right)  ,$ there exists $\left(  x_{1},x_{2},\dots,x_{n}\right)  \in T$ such
that $F=\{\max\,$supp\emph{\thinspace}$x_{i}\}_{i=1}^{n}.$ Since $\left(
x_{1},x_{2},\dots,x_{n}\right)  \in T,$ $\left(  x_{1},x_{2},\dots
,x_{n}\right)  \overset{H}{\succeq}\ell^{1}\left(  \left|  F\right|  \right)
$-basis$.$ Thus%
\[
K\rho\left(  \sum_{k=1}^{n}a_{k}e_{i_{k}}\right)  +4\varepsilon_{0}\sum
_{k=1}^{n}\left|  a_{k}\right|  \geq\frac{1}{H}\sum_{k=1}^{n}\left|
a_{k}\right|  .
\]
Hence%
\[
\rho\left(  \sum_{k=1}^{n}a_{k}e_{i_{k}}\right)  \geq\frac{1}{2KH}\sum
_{k=1}^{n}\left|  a_{k}\right|  .
\]
Since it is clear that $(e_{k})_{k\in G}\overset{1}{\succeq}(e_{k})_{k\in F}$
whenever $G$ is a spreading of $F$,
it follows that
\begin{equation}
\label{E6}\rho\left(  \sum_{k\in G}a_{k}e_{k}\right)  \geq\frac{1}{2KH}%
\sum_{k\in G}\left|  a_{k}\right|
\end{equation}
for all $G\in\mathcal{G}\left(  T\right)  .$ Assume that $\gamma_{n}\neq0$ for
some $n.$ Choose $m\in\mathbb{N}$ such that $\pi_{m}<1/(4KH)$ and $\gamma
_{m}\neq0.$ If $\mathcal{G}\left(  T\right)  $ is compact, then $\iota\left(
\mathcal{G}\left(  T\right)  \right)  >\omega^{\gamma+\gamma_{m}\cdot\omega}$
by Proposition \ref{L8}. Since $\mathcal{G}\left(  T\right)  $ is regular, the
same holds for $\mathcal{G}\left(  T\right)  \cap\left[  L\right]  ^{<\infty}$
for any $L\in\left[  \mathbb{N}\right]  .$ Thus by \cite[Corollary 1.2]{G},
there exists $L\in\left[  \mathbb{N}\right]  $ such that $\mathcal{S}%
_{\gamma+\gamma_{m}\cdot\omega}\cap\left[  L\right]  ^{<\infty}\subseteq
\mathcal{G}\left(  T\right)  .$ The same conclusion clearly holds if
$\mathcal{G}(T)$ is not compact. Hence inequality (\ref{E6}) holds for all
$\left(  a_{k}\right)  \in c_{00}$ and all $G\in\mathcal{S}_{\gamma+\gamma
_{m}\cdot\omega}\cap\left[  L\right]  ^{<\infty}.$
Now, defining $\mathcal{M}$ to be as in Lemma \ref{P3} corresponding to $m$,
\[
\iota\left(  \mathcal{M}\left[  \mathcal{F}_{0}\cup\mathcal{G}_{\varepsilon
_{0}}\right]  \right)  \leq\iota\left(  \mathcal{F}_{0}\cup\mathcal{G}%
_{\varepsilon_{0}}\right)  \cdot\iota\left(  \mathcal{M}\right)  =\left(
\alpha_{0}\vee\iota\left(  \mathcal{G}_{\varepsilon_{0}}\right)  \right)
\cdot\alpha_{m}^{m}<\omega^{\gamma+\gamma_{m}\cdot m+1}.
\]
Using \cite[Corollary 1.2]{G} again, we obtain $M\in\left[  L\right]  $ such
that $\mathcal{M}\left[  \mathcal{F}_{0}\cup\mathcal{G}_{\varepsilon_{0}%
}\right]  \cap\left[  M\right]  ^{<\infty}\subseteq\mathcal{S}_{\gamma
+\gamma_{m}\cdot m+1}.$ It follows from \cite[Proposition 3.6]{OTW} that there
are $F\in\mathcal{S}_{\gamma+\gamma_{m}\cdot\omega}\left(  M\right)  $ and
$\left(  a_{j}\right)  _{j\in F}\subseteq\mathbb{R}^{+}$ such that $\sum_{j\in
F}a_{j}=1$ and if $G\subseteq F$ with $G\in\mathcal{S}_{\gamma+\gamma_{m}\cdot
m+1},$ then $\sum_{j\in G}a_{j}<\frac{1}{4p\left(  m\right)  KH}.$ Note that
$F\in\mathcal{S}_{\gamma+\gamma_{m}\cdot\omega}\cap\left[  M\right]
^{<\infty} \subseteq\mathcal{G}(T).$ Consider $x=\sum_{j\in F}a_{j}e_{j}.$ By
Lemma \ref{P3},%
\begin{align*}
\rho\left(  x\right)   &  \leq\pi_{m}\left\|  x\right\|  _{\ell^{1}}+p\left(
m\right)  \left\|  x\right\|  _{\mathcal{M}\left[  \mathcal{F}_{0}%
\cup\mathcal{G}_{\varepsilon_{0}}\right]  }\\
&  \leq\pi_{m}+p\left(  m\right)  \left\|  x\right\|  _{\mathcal{S}%
_{\gamma+\gamma_{m}\cdot m+1}}\\
&  <\frac{1}{4KH}+\frac{1}{4KH}=\frac{1}{2KH},
\end{align*}
contrary to (\ref{E6}). This proves the proposition in case $\gamma_{n}\neq0$
for some $n.$

If $\gamma_{n}=0$ for all $n,$ then $\alpha_{n}^{\omega}=\omega$ for all $n.$
(Recall that we assume $\alpha_{n}>1$ for all $n\in\mathbb{N}$.) Write
$\alpha_{0}\vee\iota\left(  \mathcal{G}_{\varepsilon_{0}}\right)
=\omega^{\lambda_{1}}\cdot m_{1}+\dots+\omega^{\lambda_{k}}\cdot m_{k}$ in
Cantor normal form. Then
\[
\iota\left(  \mathcal{G}\left(  T\right)  \right)  \geq o\left(  T\right)
>\left[  \alpha_{0}\vee\iota\left(  \mathcal{G}_{\varepsilon_{0}}\right)
\right]  \cdot\omega=\omega^{\lambda_{1}+1}.
\]
By \cite[Corollary 1.2]{G}, there exists $L\in\left[  \mathbb{N}\right]  $
such that $\mathcal{S}_{\lambda_{1}+1}\cap\left[  L\right]  ^{<\infty
}\subseteq\mathcal{G}\left(  T\right)  .$ Hence, for all $\left(
a_{k}\right)  \in c_{00}$ and all $G\in\mathcal{S}_{\lambda_{1}+1}\cap\left[
L\right]  ^{<\infty},$ inequality (\ref{E6}) holds.
Choose $m\in\mathbb{N}$ such that $\pi_{m}<1/(4KH)$ and define $\mathcal{M}$
as in Lemma \ref{P3} corresponding to $m$. Then
\[
\iota\left(  \mathcal{M}\left[  \mathcal{F}_{0}\cup\mathcal{G}_{\varepsilon
_{0}}\right]  \right)  =\left[  \alpha_{0}\vee\iota\left(  \mathcal{G}%
_{\varepsilon_{0}}\right)  \right]  \cdot r<\omega^{\lambda_{1}}\cdot\left(
m_{1}+1\right)  r
\]
for some $r\in\mathbb{N}.$ Applying \cite[Theorem 1.1]{G}, there exists
$M\in\left[  L\right]  $ such that $\mathcal{M}\left[  \mathcal{F}_{0}%
\cup\mathcal{G}_{\varepsilon_{0}}\right]  \cap\left[  M\right]  ^{<\infty
}\subseteq\left(  \mathcal{S}_{\lambda_{1}}\right)  ^{\left(  m_{1}+1\right)
r}.$ By \cite[Proposition 3.6]{OTW}, there exist $F\in\mathcal{S}_{\lambda
_{1}+1}\left(  M\right)  \subseteq\mathcal{S}_{\lambda_{1}+1}\cap\left[
M\right]  ^{<\infty}\subseteq\mathcal{S}_{\lambda_{1}+1}\cap\left[  L\right]
^{<\infty}$ and $\left(  a_{j}\right)  _{j\in F}\subseteq\mathbb{R}^{+}$ such
that $\sum_{j\in F}a_{j}=1$ and if $G\subseteq F$ with $G\in\mathcal{S}%
_{\lambda_{1}},$ then $\sum_{j\in G}a_{j}<\frac{1}{4p(m)KH(m_{1}+1)r}.$
Consider $x=\sum_{j\in F}a_{j}e_{j}.$ By Lemma \ref{P3},%
\begin{align*}
\rho\left(  x\right)   &  \leq\pi_{m}\left\|  x\right\|  _{\ell^{1}}+p\left(
m\right)  \left\|  x\right\|  _{\mathcal{M}\left[  \mathcal{F}_{0}%
\cup\mathcal{G}_{\varepsilon_{0}}\right]  }\\
&  \leq\pi_{m}+p\left(  m\right)  \left\|  x\right\|  _{\left(  \mathcal{S}%
_{\lambda_{1}}\right)  ^{\left(  m_{1}+1\right)  r}}\\
&  \leq\pi_{m}+p\left(  m\right)  \left(  m_{1}+1\right)  r\left\|  x\right\|
_{\mathcal{S}_{\lambda_{1}}}\\
&  <\frac{1}{4KH}+\frac{1}{4KH}=\frac{1}{2KH},
\end{align*}
contradicting (\ref{E6}).
\end{proof}

\begin{theorem}
\label{P7}

\begin{enumerate}
\item $\alpha_{0}\cdot\sup\limits_{n\in\mathbb{N}}\alpha_{n}^{\omega}\leq
I_{b}\left(  X\right)  \leq\left(  \alpha_{0}\vee\alpha\right)  \cdot
\sup\limits_{n\in\mathbb{N}}\alpha_{n}^{\omega}.$

\item If $\alpha_{0}\geq\alpha,$ then $I_{b}\left(  X\right)  =\alpha_{0}%
\cdot\sup\limits_{n\in\mathbb{N}}\alpha_{n}^{\omega}.$

\item If $\alpha_{0}<\alpha$ and $\alpha=\alpha_{n}$ for some $n\in
\mathbb{N},$ then $I_{b}\left(  X\right)  =\alpha^{\omega}.$

\item If $\alpha_{n}<\alpha$ for all $n\in\mathbb{N}\cup\left\{  0\right\}  $
and $\alpha$ is not of the form $\alpha=\omega^{\omega^{\xi}},$ $\xi
<\omega_{1},$ then $I_{b}\left(  X\right)  =\alpha^{\omega}.$
\end{enumerate}
\end{theorem}

\begin{proof}
1.\ The first inequality follows from Proposition \ref{P8}. Since
$\mathcal{S}_{0}\subseteq\mathcal{F}_{n}$ for all $n,$ $\mathcal{F}%
_{m}\subseteq\mathcal{F}\subseteq\lbrack\mathcal{F},\mathcal{F}_{n_{1}}%
,\dots,\mathcal{F}_{n_{s}}]$ for all $m,n_{1},\dots,n_{s}\in\mathbb{N}.$ The
second inequality follows from Proposition \ref{P11} upon taking
$\mathcal{G}_{\varepsilon}=\mathcal{F}$.

2.\ and 3.\ are clear.

4.\ In this case, it is readily verified that $\alpha\cdot\sup\limits_{n\in
\mathbb{N}}\alpha_{n}^{\omega}=\sup\limits_{n\in\mathbb{N}}\alpha_{n}^{\omega
}=\alpha^{\omega}.$ The conclusion follows from 1.
\end{proof}

The following corollary answers Question 1 in \cite{JO}.

\begin{corollary}
\label{solution} If $\eta$ is a countable ordinal not of the form $\omega
^{\xi}$ for some limit ordinal $\xi<\omega_{1},$ then there exists a Banach
space $Y$ such that $I\left(  Y\right)  =\omega^{\eta}.$
\end{corollary}

\begin{proof}
Write $\eta=\omega^{\gamma_{1}}\cdot m_{1}+\dots+\omega^{\gamma_{k}}\cdot
m_{k}$ in Cantor normal form. If $\gamma_{k}$ is $0$ or a successor ordinal,
then the result follows immediately from \cite[Corollary 14]{LT}. If
$\gamma_{k}$ is a limit ordinal, let $\left(  \beta_{n}\right)  $ be a
sequence of ordinals increasing to $\gamma_{k}.$ Choose regular families
$\left(  \mathcal{F}_{n}\right)  _{n=0}^{\infty}$ such that $\alpha_{n}%
=\iota\left(  \mathcal{F}_{n}\right)  =\omega^{\omega^{\beta_{n}}}$,
$n\in\mathbb{N},$ and $\alpha_{0}=\iota\left(  \mathcal{F}_{0}\right)
=\omega^{\omega^{\gamma_{1}}\cdot m_{1}+\dots+\omega^{\gamma_{k}}\cdot
(m_{k}-1)}.$ Then $\alpha=\sup\limits_{n\in\mathbb{N}}\alpha_{n}%
=\sup\limits_{n\in\mathbb{N}}\omega^{\omega^{\beta_{n}}}=\omega^{\omega
^{\gamma_{k}}}\leq\alpha_{0}$ as $k>1$ or $m_{k}>1.$ Let $Y=T\left(
\mathcal{F}_{0},\left(  \theta_{n},\mathcal{F}_{n}\right)  _{n=1}^{\infty
}\right)  .$ By 2.\ in Theorem \ref{P7}, $I_{b}\left(  Y\right)  =\alpha
_{0}\cdot\sup\limits_{n\in\mathbb{N}}\alpha_{n}^{\omega}=\omega^{\omega
^{\gamma_{1}}\cdot m_{1}+\dots+\omega^{\gamma_{k}}\cdot m_{k}}=\omega^{\eta}.$
Finally, since $I_{b}\left(  Y\right)  \geq\omega^{\omega},$ $I\left(
Y\right)  =I_{b}\left(  Y\right)  =\omega^{\eta}$ by \cite[Corollary 5.13]{JO}.
\end{proof}

\section{\label{S4}Attaining the upper bound}

Henceforth, we shall consider only the case where $\alpha_{n}<\alpha$ for all
$n\in\mathbb{N}\cup\left\{  0\right\}  $ and $\alpha$ is of the form
$\omega^{\omega^{\xi}}.$ Under these conditions, Theorem \ref{P7} yields the
estimate
\[
\omega^{\omega^{\xi}}\leq I_{b}\left(  X\right)  \leq\omega^{\omega^{\xi}%
\cdot2}.
\]
The next theorem gives a sufficient condition for the upper estimate to be
attained. Given $m\in\mathbb{N}$ and $\varepsilon>0,$ define
\begin{multline*}
\gamma=\gamma\left(  \varepsilon,m\right)  =\max\{ \ell( \alpha_{0}\cdot
\alpha_{n_{s}}\dots\alpha_{n_{1}}) :\\
\varepsilon\theta_{n_{1}}\theta_{n_{2}}\dots\theta_{n_{s}}>\theta_{m}\} \text{
($\max\emptyset=0$).}%
\end{multline*}

\begin{theorem}
\label{T2}Assume $\xi\neq0.$ If there exists $\varepsilon>0$ such that\ for
all $\beta<\omega^{\xi},$ there exists $m\in\mathbb{N}$ satisfying
$\gamma\left(  \varepsilon,m\right)  +2+\beta<\ell\left(  \alpha_{m}\right)
$, then $I_{b}\left(  X\right)  =\omega^{\omega^{\xi}\cdot2}$.
\end{theorem}

Before giving the proof of Theorem \ref{T2}, let us observe an interesting corollary.

\begin{corollary}
If $\xi$ is a limit ordinal, then $I_{b}\left(  X\right)  =\omega^{\omega
^{\xi}\cdot2}.$
\end{corollary}

\begin{proof}
Since $\xi$ is a limit ordinal, the sequence $\left(  \ell\left(  \ell\left(
\alpha_{n}\right)  \right)  \right)  $ converges to $\xi.$ Hence for all
$\beta<\xi,$ there exists $m\in\mathbb{N}$ such that $\ell\left(  \ell\left(
\alpha_{m}\right)  \right)  >\beta\vee\ell\left(  \ell\left(  \alpha
_{m-1}\right)  \right)  \vee\ell\left(  \ell\left(  \alpha_{0} \right)
\right)  .$ Suppose $\theta_{n_{1}}\dots\theta_{n_{s}}>\theta_{m}$ for some
$n_{1},\dots,n_{s}\in\mathbb{N}.$ Then $n_{1},\dots,n_{s}<m.$ Now for all
$1\leq i\leq s,$ $\alpha_{n_{i}}\leq\omega^{\ell\left(  \alpha_{n_{i}}\right)
+1}.$ Thus%
\[
\alpha_{0}\cdot\alpha_{n_{s}}\dots\alpha_{n_{1}}\leq\omega^{\ell\left(
\alpha_{0}\right)  +1+\ell\left(  \alpha_{n_{s}}\right)  +1+\dots+\ell\left(
\alpha_{n_{1}}\right)  +1}.
\]
Therefore%
\begin{align*}
\ell(\alpha_{0}  &  \cdot\alpha_{n_{s}}\dots\alpha_{n_{1}}) +2+\omega^{\beta
}\\
&  \leq\ell(\alpha_{0})+1+\ell(\alpha_{n_{s}})+1+\dots+\ell(\alpha_{n_{1}})
+1+2+\omega^{\beta}\\
&  <\omega^{\ell(\ell(\alpha_{m}))}\leq\ell(\alpha_{m}) .
\end{align*}
Applying Theorem \ref{T2} with $\varepsilon=1$ yields the required result.
\end{proof}

\begin{lemma}
\label{LB1}Let $m\in\mathbb{N}$ and $\varepsilon>0$ be given$.$ Then for all
$M\in\left[  \mathbb{N}\right]  ,$ there exists $x\in c_{00}$ satisfying
$\left\|  x\right\|  \leq1+\frac{1}{\varepsilon},$ $\left\|  x\right\|
_{\ell^{1}}=\frac{1}{\theta_{m}},$ and $\emph{supp}$ $x\in\mathcal{S}%
_{\gamma+2}\cap\left[  M\right]  ^{<\infty},$ where $\gamma=\gamma\left(
\varepsilon,m\right)  $ is as defined above$.$
\end{lemma}

\begin{proof}
Let $\mathcal{N}=\left\{  \left(  n_{1},\dots,n_{s}\right)  :\varepsilon
\theta_{n_{1}}\dots\theta_{n_{s}}>\theta_{m}\right\}  .$ Clearly $\mathcal{N}$
is a finite set. Denote its cardinality by $c.$ By assumption, there exists
$L\in\left[  M\right]  ^{<\infty}$ such that $\left[  \mathcal{F}_{n_{1}%
},\dots,\mathcal{F}_{n_{s}},\mathcal{F}_{0}\right]  \cap\left[  L\right]
^{<\infty}\subseteq\mathcal{S}_{\gamma+1}$ for all $\left(  n_{1},\dots
,n_{s}\right)  \in\mathcal{N}$ (cf. \cite{G})$.$ By \cite[Proposition
3.6]{OTW}, there exists $y\in c_{00},$ $\left\|  y\right\|  _{\ell^{1}}=1$
such that supp\thinspace$y\in\mathcal{S}_{\gamma+2}\cap\left[  L\right]
^{<\infty}$ and $\left|  \left|  y\right|  \right|  _{\mathcal{S}_{\gamma+1}%
}\leq\theta_{m}/c.$ Let $x=y/\theta_{m}.$ Then $\left|  \left|  x\right|
\right|  _{\ell^{1}}=\frac{1}{\theta_{m}}$ and supp $x\in\mathcal{S}%
_{\gamma+2}\cap\left[  M\right]  ^{<\infty}.$ Choose a complete $\left(
\mathcal{F}_{n}\right)  $-admissible tree $\mathcal{T}$ such that $\left|
\left|  x\right|  \right|  =\mathcal{T}x.$ Denote by $\mathcal{L}\left(
\mathcal{T}\right)  $ the set of all leaves of $\mathcal{T}.$ For a fixed
$\left(  n_{1},\dots,n_{s}\right)  \in\mathcal{N},$ the set $\left\{
E\in\mathcal{L}\left(  \mathcal{T}\right)  :h\left(  E\right)  =\left(
0,n_{1},\dots,n_{s}\right)  \right\}  $ is $\left[  \mathcal{F}_{n_{1}}%
,\dots,\mathcal{F}_{n_{s}}\right]  $-admissible. Since supp\thinspace
$x\in\left[  L\right]  ^{<\infty},$ we conclude by the choice of $L$ that
\[
\sum_{\substack{E\in\mathcal{L}\left(  \mathcal{T}\right)  \\\,h\left(
E\right)  =\left(  0,n_{1},\dots,n_{s}\right)  }}\left|  \left|  Ex\right|
\right|  _{\mathcal{F}_{0}}\leq\left\|  x\right\|  _{\left[  \mathcal{F}%
_{n_{1}},\dots,\mathcal{F}_{n_{s}},\mathcal{F}_{0}\right]  }\leq\left\|
x\right\|  _{\mathcal{S}_{\gamma+1}}.
\]
Therefore
\begin{align*}
\left|  \left|  x\right|  \right|   &  \leq\sum_{\substack{E\in\mathcal{L}%
\left(  \mathcal{T}\right)  \\\varepsilon t\left(  E\right)  \leq\theta_{m}%
\ }}t\left(  E\right)  \left|  \left|  Ex\right|  \right|  _{\mathcal{F}_{0}%
}+\sum_{\substack{E\in\mathcal{L}\left(  \mathcal{T}\right)  \\\varepsilon
t\left(  E\right)  >\theta_{m}\ }}t\left(  E\right)  \left|  \left|
Ex\right|  \right|  _{\mathcal{F}_{0}}\\
&  \leq\frac{\theta_{m}}{\varepsilon}\left|  \left|  x\right|  \right|
_{\ell^{1}}+\sum_{\left(  n_{1},\dots,n_{s}\right)  \in\mathcal{N}}%
\theta_{n_{1}}\theta_{n_{2}}\dots\theta_{n_{j}}\sum_{\substack{E\in
\mathcal{L}\left(  \mathcal{T}\right)  \\h\left(  E\right)  =\left(
0,n_{1},\dots,n_{s}\right)  }}\left|  \left|  Ex\right|  \right|
_{\mathcal{F}_{0}}\\
&  \leq\frac{1}{\varepsilon}+\sum_{\left(  n_{1},\dots,n_{s}\right)
\in\mathcal{N}}\left|  \left|  x\right|  \right|  _{\mathcal{S}_{\gamma+1}%
}\leq\frac{1}{\varepsilon}+\frac{c}{\theta_{m}}\left|  \left|  y\right|
\right|  _{\mathcal{S}_{\gamma+1}}\leq1+\frac{1}{\varepsilon}.
\end{align*}
\end{proof}

\begin{lemma}
\label{LB2}Under the assumptions of Theorem \ref{T2}, there exists a strictly
increasing sequence $\left(  q_{k}\right)  \subseteq\mathbb{N}$ such that for
all $F\in\mathcal{S}_{\omega^{\xi}},$ there are normalized vectors $\left(
x_{k}\right)  _{k\in F}$ with \emph{supp} $x_{k}\subseteq\lbrack q_{k}%
,q_{k+1})$ for all $k\in F$ and
\[
\left|  \left|  \sum_{k\in F}a_{k}x_{k}\right|  \right|  \geq\frac
{\varepsilon}{1+\varepsilon}\sum_{k\in F}\left|  a_{k}\right|
\]
for all $\left(  a_{k}\right)  \in c_{00}.$
\end{lemma}

\begin{proof}
Since $\xi\neq0,$ $\omega^{\xi}$ is a limit ordinal. Suppose that
$\mathcal{S}_{\omega^{\xi}}$ is defined by the sequence $\left(  \beta
_{k}\right)  $ increasing to $\omega^{\xi}$. For each $k,$ choose $m_{k}%
\in\mathbb{N}$ such that $\gamma\left(  \varepsilon,m_{k}\right)  +2+\beta
_{k}<\ell\left(  \alpha_{m_{k}}\right)  .$ Write $\gamma_{k}=\gamma\left(
\varepsilon,m_{k}\right)  .$ Using Lemma \ref{LB1}, obtain a strictly
increasing sequence $\left(  q_{k}\right)  _{k=1}^{\infty}$ in $\mathbb{N}$
and $\left.  \left(  x_{k}^{i}\right)  _{k=1}^{i}\right.  _{i=1}^{\infty
}\subseteq c_{00}$ such that $\left\|  x_{k}^{i}\right\|  _{\ell^{1}}=\frac
{1}{\theta_{m_{k}}},$ $\left\|  x_{k}^{i}\right\|  \leq1+\frac{1}{\varepsilon
},$ supp $x_{k}^{i}\subseteq\lbrack q_{i},q_{i+1}),$ and supp $x_{k}^{i}%
\in\mathcal{S}_{\gamma_{k}+2}\cap[M_{i}]^{<\infty},$ where \thinspace$M_{i}%
\in\left[  \mathbb{N}\right]  $ is chosen so that $M_{i+1}\subseteq M_{i}%
\cap\lbrack q_{i+1},\infty)$ and
\[
\bigcup_{j=1}^{i}\mathcal{S}_{\beta_{j}}[\mathcal{S}_{\gamma_{i}+2}]
\cap[M_{i}] ^{<\infty}\subseteq\mathcal{F}_{m_{i}}.
\]
Note that this choice is possible by \cite{G} since
\[
\iota( \cup_{j=1}^{i}\mathcal{S}_{\beta_{j}}\left[  \mathcal{S}_{\gamma_{i}%
+2}\right]  ) =\omega^{\gamma_{i}+2+\beta_{i}}<\omega^{\ell( \alpha_{m_{i}})
}\leq\alpha_{m_{i}}=\iota\left(  \mathcal{F}_{m_{i}}\right)  .
\]
If $F=\left\{  i_{1},\dots,i_{r}\right\}  \in\mathcal{S}_{\omega^{\xi}}$,
$i_{1} < \dots< i_{r}$, then $F\in\mathcal{S}_{\beta_{k}}$ for some $k\leq
i_{1}.$ Consider the block basic sequence $(x_{i_{1}}^{i_{1}},x_{i_{1}}%
^{i_{2}},\dots,x_{i_{1}}^{i_{r}}).$ By choice, supp $x_{i_{1}}^{i_{j}}%
\in\mathcal{S}_{\gamma_{i_{1}}+2}\cap[M_{i_{j}}]^{<\infty}$ and supp
$x_{i_{1}}^{i_{j}}\subseteq\lbrack q_{i_{j}},q_{i_{j}+1}),$ $1\leq j\leq r.$
Moreover, the set $\{q_{i_{1}},\dots,q_{i_{r}}\} $ is a spreading of
$\{i_{1},\dots,i_{r}\} = F$ and hence belongs to $\mathcal{S}_{\beta_{k}}.$
Thus
\[
\bigcup_{j=1}^{r}\text{supp\thinspace}x_{i_{1}}^{i_{j}}\in\mathcal{S}%
_{\beta_{k}}[ \mathcal{S}_{\gamma_{i_{1}}+2}] \cap[ M_{i_{1}}] ^{<\infty
}\subseteq\mathcal{F}_{m_{i_{1}}}.
\]
Therefore, given any $\left(  a_{j}\right)  \in c_{00},$
\begin{align*}
\| \sum_{j=1}^{r}a_{j}x_{i_{1}}^{i_{j}}\|  &  \geq\theta_{m_{i_{1}}}\|
\sum_{j=1}^{r}a_{j}x_{i_{1}}^{i_{j}}\| _{\ell^{1}}\\
&  =\theta_{m_{i_{1}}}\sum_{j=1}^{r}| a_{j}| \| x_{i_{1}}^{i_{j}}\| _{\ell
^{1}}\\
&  =\theta_{m_{i_{1}}}\sum_{j=1}^{r}\left|  a_{j}\right|  \frac{1}%
{\theta_{m_{i_{1}}}}=\sum_{j=1}^{r}\left|  a_{j}\right|  .
\end{align*}
Normalizing the sequence $( x_{i_{1}}^{i_{1}},x_{i_{1}}^{i_{2}},\dots
,x_{i_{1}}^{i_{r}}) $ yields the desired result.
\end{proof}

\begin{lemma}
\label{LB3}Suppose the assumptions of Theorem \ref{T2} hold. Then there exists
$\left(  q_{k}\right)  \subseteq\mathbb{N}$ such that whenever $F\in
\mathcal{F}_{\alpha_{n}}\left[  \mathcal{S}_{\omega^{\xi}}\right]  $ for some
$n\in\mathbb{N},$ there are normalized vectors $\left(  x_{k}\right)  _{k\in
F},$ \emph{supp} $x_{k}\subseteq\lbrack q_{k},q_{k+1}),$ satisfying
\[
\left|  \left|  \sum_{k\in F}a_{k}x_{k}\right|  \right|  \geq\frac
{\varepsilon\theta_{n}}{1+\varepsilon}\sum_{k\in F}\left|  a_{k}\right|
\]
for all $\left(  a_{k}\right)  \in c_{00}.$
\end{lemma}

\begin{proof}
Choose $\left(  q_{k}\right)  $ using Lemma \ref{LB2}. If $F\in\mathcal{F}%
_{\alpha_{n}}\left[  \mathcal{S}_{\omega^{\xi}}\right]  $ for some
$n\in\mathbb{N},$ write $F=\bigcup_{j=1}^{s}F_{j},$ with $F_{1}<\dots<F_{s},$
$F_{j}\in\mathcal{S}_{\omega^{\xi}}$, $1\leq j\leq s,$ and $\left\{  \min
F_{j}\right\}  _{j=1}^{s}\in\mathcal{F}_{\alpha_{n}}.$ For all $1\leq j\leq
s,$ there exist normalized vectors $\left(  x_{k}\right)  _{k\in F_{j}}$ such
that supp $x_{k}\subseteq\lbrack q_{k},q_{k+1})$ for all $k\in F_{j}$ and
$\left|  \left|  \sum_{k\in F_{j}}a_{k}x_{k}\right|  \right|  \geq
\frac{\varepsilon}{1+\varepsilon}\sum_{k\in F_{j}}\left|  a_{k}\right|  $ for
any $\left(  a_{k}\right)  \in c_{00}.$ Therefore,
\begin{align*}
\left\|  \sum_{k\in F}a_{k}x_{k}\right\|   &  =\left\|  \sum_{j=1}^{s}\left(
\sum_{k\in F_{j}}a_{k}x_{k}\right)  \right\| \\
&  \geq\theta_{n}\sum_{j=1}^{s}\left\|  E_{j}\sum_{j=1}^{s}\left(  \sum_{k\in
F_{j}}a_{k}x_{k}\right)  \right\|  ,\text{ where }E_{j}=\bigcup_{k\in F_{j}%
}\text{supp\thinspace}x_{k}\\
&  =\theta_{n}\sum_{j=1}^{s}\left\|  \sum_{k\in F_{j}}a_{k}x_{k}\right\| \\
&  \geq\frac{\varepsilon\theta_{n}}{1+\varepsilon}\sum_{k\in F}\left|
a_{k}\right|
\end{align*}
for any $\left(  a_{k}\right)  \in c_{00}.$
\end{proof}

To complete the proof of Theorem \ref{T2}, we apply a compactness argument to
condense the block basic sequences obtained in Lemma \ref{LB3} into a tree.
Let $Y$ be a set and let $\left(  A_{n}\right)  _{n=1}^{\infty}$ be a sequence
of pairwise disjoint finite subsets of $Y.$ Suppose that a given set
\[
\mathcal{X}\subseteq\bigcup_{\emptyset\neq F\in\left[  \mathbb{N}\right]
^{<\infty}}\left(  \prod_{n\in F}A_{n}\right)
\]
is hereditary in the sense that $\left(  x_{n}\right)  _{n\in G}\in
\mathcal{X}$ whenever $\left(  x_{n}\right)  _{n\in F}\in\mathcal{X}$ and
$\emptyset\neq G\subseteq F.$

\begin{proposition}
\label{P20}Let $\mathcal{H}\subseteq\left[  \mathbb{N}\right]  ^{<\infty}$ be
a regular family with $\omega_{1}>\iota\left(  \mathcal{H}\right)  \geq
\alpha\geq1.$ Suppose for all nonempty $F\in\mathcal{H},$ there exists
$\left(  x_{n}\right)  _{n\in F}\in\mathcal{X}.$ Then there exists a tree $T$
on $Y$ such that $T\subseteq\mathcal{X}$ and $o\left(  T\right)  \geq\alpha.$
\end{proposition}

\begin{proof}
Assume that $\mathcal{H}$ is regular and nonempty. There exists $n_{0}%
\in\mathbb{N}$ such that $\left\{  n\right\}  \in\mathcal{H}$ for all $n\geq
n_{0}.$ By hypothesis, there exists $\left(  x_{n}\right)  \in\mathcal{X}$ for
all $n\geq n_{0}.$ Let $T=\left\{  \left(  x_{n}\right)  :n\geq n_{0}\right\}
.$ Then $T\subseteq\mathcal{X}$ and $o\left(  T\right)  \geq1.$

Suppose the proposition is true for some $\alpha\geq1.$ Let $\mathcal{H}%
\subseteq\left[  \mathbb{N}\right]  ^{<\infty}$ be a regular family satisfying
the hypothesis such that $\omega_{1}>\iota\left(  \mathcal{H}\right)
\geq\alpha+1.$ Pick a singleton set $\{n_{0}\}\in\mathcal{H}^{\left(
\alpha\right)  }$ and let%
\[
\mathcal{G}=\left\{  G\in\left[  \mathbb{N}\right]  ^{<\infty}: n_{0}<G,
\{n_{0}\}\cup G\in\mathcal{H}\right\}  .
\]
Then $\mathcal{G}$ is regular and $\iota\left(  \mathcal{G}\right)  \geq
\alpha\geq1.$ Correspondingly, let
\[
\begin{split}
\mathcal{Y} =\{ ( x_{n}) _{n\in G}:\emptyset\neq\  &  G\in\mathcal{G}\text{,
there exists } (x_{n_{0}})\\
&  \text{such that }( x_{n}) _{n\in\{n_{0}\}\cup G}\in\mathcal{X}\} .
\end{split}
\]
Since $\mathcal{X}$ is hereditary, so is $\mathcal{Y}$. Let a nonempty set
$G\in\mathcal{G}$ be given. Then there exists $\left(  x_{n}\right)
_{n\in\{n_{0}\}\cup G}\in\mathcal{X}$ such that $\left(  x_{n}\right)  _{n\in
G}\in\mathcal{Y}.$ By the inductive hypothesis, there exists a tree $T_{0}$ on
$Y$ such that $T_{0}\subseteq\mathcal{Y}$ and $o\left(  T_{0}\right)
\geq\alpha.$
List the elements in
$A_{n_{0}}$ as $( z_{n_{0}}^{1}), \dots, ( z_{n_{0}}^{p}).$ Let $M$ be the
collection of maximal nodes of $T_{0}.$ If $\left(  x_{n}\right)  _{n\in G}\in
M,$ there exists $i$, $1\leq i\leq p,$ such that $( z_{n_{0}}^{i})\cup
(x_{n})_{n\in G}\in\mathcal{X}.$ Partition $M$ into $\bigcup_{i=1}^{p}M_{i}$
so that $\left(  x_{n}\right)  _{n\in G}\in M_{i}$ implies $( z_{n_{0}}^{i})
\cup\left(  x_{n}\right)  _{n\in G}\in\mathcal{X}.$ Now let
$T_{i}$ be the subtree of $T_{0}$ consisting of all nodes in $M_{i}$ and their
ancestors.
By \cite[Lemma 5.10]{JO}, there exists $i$ such that $o\left(  T_{i}\right)
\geq\alpha.$ Define
\[
T=\left\{  (z_{n_{0}}^{i})\cup\left(  x_{n}\right)  _{n\in H}:\left(
x_{n}\right)  _{n\in H}\in T_{i}\right\}  .
\]
Then $T$ is a tree on $Y$ such that $T\subseteq\mathcal{X}$ and $o\left(
T\right)  \geq\alpha+1.$

Suppose $\alpha$ is a countable limit ordinal and the result holds for all
$1\leq\beta<\alpha.$ Let $\mathcal{H}\subseteq\left[  \mathbb{N}\right]
^{<\infty}$ be a regular family of finite subsets of $\mathbb{N}$ satisfying
the hypothesis such that $\iota\left(  \mathcal{H}\right)  \geq\alpha.$ If
$1\leq\beta<\alpha,$ then $\iota\left(  \mathcal{H}\right)  \geq\beta\geq1.$
Hence there exists a tree $T_{\beta}$ on $Y$ such that $T_{\beta}%
\subseteq\mathcal{X}$ and $o\left(  T_{\beta}\right)  \geq\beta.$ Clearly the
tree $T=\bigcup_{\beta<\alpha}T_{\beta}$ satisfies the requirements of the proposition.
\end{proof}

\begin{proof}
[Proof of Theorem \ref{T2}]In view of 1.\ in Theorem \ref{P7}, it suffices to
show that $I_{b}\left(  X\right)  \geq\omega^{\omega^{\xi}}\cdot\alpha_{n}$
for all $n\in\mathbb{N}.$ In order to set up to apply Proposition \ref{P20},
let $Y=X.$ Choose a sequence $\left(  q_{k}\right)  $ as in Lemma \ref{LB3}
and fix $n\in\mathbb{N.}$ Let $A_{k}$ be a finite $\frac{\varepsilon\theta
_{n}}{2\left(  1+\varepsilon\right)  }$-net of the unit sphere of $\left[
e_{j}\right]  _{j=q_{k}}^{q_{k+1}-1}$ for each $k\in\mathbb{N}.$ Define
$c_{n}=\frac{\varepsilon\theta_{n}}{2\left(  1+\varepsilon\right)  }$ and set%
\[
\mathcal{X}=\{ \left(  y_{k}\right)  _{F}: \emptyset\neq F\in\mathcal{F}%
_{\alpha_{n}}\left[  \mathcal{S}_{\omega^{\xi}}\right]  ,\text{ }y_{k}\in
A_{k},\,\text{ }\left(  y_{k}\right)  \overset{c_{n}}{\succeq}\ell^{1}\left(
\left|  F\right|  \right)  \!\text{-basis}\} .
\]
Clearly $\mathcal{X}$ is hereditary. According to Lemma \ref{LB3}, whenever
$F\in\mathcal{F}_{\alpha_{n}}\left[  \mathcal{S}_{\omega^{\xi}}\right]  ,$
there exist normalized vectors $\left(  x_{k}\right)  _{k\in F},$ supp
$x_{k}\subseteq\lbrack q_{k},q_{k+1}),$ such that
\[
\left|  \left|  \sum_{k\in F}a_{k}x_{k}\right|  \right|  \geq\frac
{\varepsilon\theta_{n}}{1+\varepsilon}\sum_{k\in F}\left|  a_{k}\right|
\]
for all $\left(  a_{k}\right)  \in c_{00}.$ Choose $\left(  y_{k}\right)
_{k\in F}$ such that $y_{k}\in A_{k}$ and $\left|  \left|  x_{k}-y_{k}\right|
\right|  \leq\frac{\varepsilon\theta_{n}}{2\left(  1+\varepsilon\right)  }$
for all $k\in F.$ For all $\left(  a_{k}\right)  \in c_{00},$%
\begin{align*}
\left|  \left|  \sum_{k\in F}a_{k}y_{k}\right|  \right|   &  \geq\left|
\left|  \sum_{k\in F}a_{k}x_{k}\right|  \right|  -\left|  \left|  \sum_{k\in
F}a_{k}\left(  x_{k}-y_{k}\right)  \right|  \right| \\
&  \geq\frac{\varepsilon\theta_{n}}{1+\varepsilon}\sum_{k\in F}\left|
a_{k}\right|  -\sum_{k\in F}\left|  a_{k}\right|  \left|  \left|  x_{k}%
-y_{k}\right|  \right| \\
&  \geq\frac{\varepsilon\theta_{n}}{1+\varepsilon}\sum_{k\in F}\left|
a_{k}\right|  -\frac{\varepsilon\theta_{n}}{2\left(  1+\varepsilon\right)
}\sum_{k\in F}\left|  a_{k}\right| \\
&  =\frac{\varepsilon\theta_{n}}{2\left(  1+\varepsilon\right)  }\sum_{k\in
F}\left|  a_{k}\right|  .
\end{align*}
Thus $\left(  y_{k}\right)  _{k\in F}\in\mathcal{X}.$ By Proposition
\ref{P20}, there exists a tree $T$\ on $X$ such that $T\subseteq\mathcal{X}$
and $o\left(  T\right)  \geq\iota\left(  \mathcal{F}_{\alpha_{n}}\left[
\mathcal{S}_{\omega^{\xi}}\right]  \right)  =\omega^{\omega^{\xi}}\cdot
\alpha_{n}.$ Since $T\subseteq\mathcal{X},$ it is an $\ell^{1}$-$c_{n}$-block
tree. Thus $I_{b}\left(  X\right)  \geq\omega^{\omega^{\xi}}\cdot\alpha_{n}.$
\end{proof}

In general, the converse of Theorem \ref{T2} is far from true, as the
following theorem shows.

\begin{theorem}
\label{T3}Suppose that $0<\xi<\omega_{1},$ $\left(  \alpha_{n}\right)
_{n=0}^{\infty}$ is a sequence of ordinals such that $\sup\limits_{n\in
\mathbb{N}\cup\left\{  0\right\}  }\alpha_{n}=\omega^{\omega^{\xi}}$
nontrivially (i.e., $\alpha_{n}<\omega^{\omega^{\xi}}$ for all $n$) and
$\left(  \theta_{n}\right)  _{n=1}^{\infty}$ is a nonincreasing null sequence
in $\left(  0,1\right)  $. Then there exists a sequence $\left(
\mathcal{F}_{n}\right)  _{n=0}^{\infty}$ of regular families of finite subsets
of $\mathbb{N}$ such that $\iota\left(  \mathcal{F}_{n}\right)  =\alpha_{n}$
for all $n\in\mathbb{N}\cup\left\{  0\right\}  $ and $I_{b}\left(  T\left(
\mathcal{F}_{0},\left(  \theta_{n},\mathcal{F}_{n}\right)  _{n=1}^{\infty
}\right)  \right)  =\omega^{\omega^{\xi}\cdot2}.$
\end{theorem}

\begin{proof}
The proof is similar to that of Theorem \ref{T2} once we have obtained
Proposition \ref{P24} below.
\end{proof}

\begin{lemma}
\label{L23}Suppose that $\omega\leq\beta<\omega_{1},$ where $\beta
=\omega^{\beta_{1}}\cdot k_{1}+\dots+\omega^{\beta_{m}}\cdot k_{m}$ in Cantor
normal form, and $g:\mathbb{N}\to\mathbb{N}$ is a function increasing to
$\infty$. There exist regular families $\mathcal{G}$ and $\mathcal{H}$ such
that $\omega\cdot\iota\left(  \mathcal{G}\right)  =\omega^{\beta_{1}}\cdot
k_{1}$, $\mathcal{S}_{0}\subseteq\mathcal{G}$ and $\iota\left(  \mathcal{H}%
\right)  =\omega^{\beta_{2}}\cdot k_{2}+\dots+\omega^{\beta_{m}}\cdot k_{m}.$
In particular, $\iota\left(  \left(  \mathcal{H},\mathcal{G}\left[
\mathcal{S}_{1}^{g}\right]  \right)  \right)  =\beta.$ (If $m = 1$, take
$\mathcal{H} = \emptyset$.)
\end{lemma}

\begin{proof}
Note that $\beta_{1}>0$ since $\beta\geq\omega.$ Define%
\[
\mathcal{G}=\left\{
\begin{array}
[c]{cc}%
\left(  \mathcal{S}_{\beta_{1}-1}\right)  ^{k_{1}} & \text{if }0<\beta
_{1}<\omega\\
\left(  \mathcal{S}_{\beta_{1}}\right)  ^{k_{1}} & \text{if }\omega\leq
\beta_{1}<\omega_{1}%
\end{array}
\right.
\]
and $\mathcal{H=}( ( \mathcal{S}_{\beta_{m}}) ^{k_{m}},\dots,( \mathcal{S}%
_{\beta_{2}}) ^{k_{2}}) .$ Clearly $\iota( \mathcal{H}) =\omega^{\beta_{2}%
}\cdot k_{2}+\dots+\omega^{\beta_{m}}\cdot k_{m}$ and
\[
\omega\cdot\iota\left(  \mathcal{G}\right)  =\left\{
\begin{array}
[c]{cc}%
\omega\cdot\omega^{\beta_{1}-1}\cdot k_{1} & \text{if }0<\beta_{1}<\omega\\
\omega\cdot\omega^{\beta_{1}}\cdot k_{1} & \text{if }\omega\leq\beta
_{1}<\omega_{1}%
\end{array}
\right.  = \omega^{\beta_{1}}\cdot k_{1}.
\]
\end{proof}

If $\beta$ is a nonzero countable ordinal whose Cantor normal form is
$\omega^{\beta_{1}}\cdot k_{1}+\dots+\omega^{\beta_{m}}\cdot k_{m}$, write
$\mathcal{R}_{\beta}$ for the family $( ( \mathcal{S}_{\beta_{m}}) ^{k_{m}%
},\dots,( \mathcal{S}_{\beta_{1}}) ^{k_{1}})$.

\begin{proposition}
\label{P24}Under the hypotheses of Theorem \ref{T3}, there exist regular
families $\left(  \mathcal{F}_{n}\right)  _{n=0}^{\infty}$\ and $\mathcal{G}$
with $\iota\left(  \mathcal{F}_{n}\right)  =\alpha_{n},$ $\iota\left(
\mathcal{G}\right)  =\omega^{\omega^{\xi}},$ and $\left(  q_{m}\right)
\subseteq\mathbb{N}$ such that for all $n\in\mathbb{N}$ and all $F\in$
$\mathcal{F}_{n}\left[  \mathcal{G}\right]  ,$ there is a normalized sequence
$\left(  x_{m}\right)  _{m\in F}$ such that \emph{supp} $x_{m}\subseteq\lbrack
q_{m},q_{m+1})$ and
\[
\left\|  \sum_{m\in F}a_{m}x_{m}\right\|  \geq\frac{\theta_{n}}{2}\sum_{m\in
F}\left|  a_{m}\right|
\]
for all $\left(  a_{m}\right)  \in c_{00}.$ Here the norm $\left\|
\cdot\right\|  $ is taken in the space $T\left(  \mathcal{F}_{0},\left(
\theta_{n},\mathcal{F}_{n}\right)  _{n=1}^{\infty}\right)  .$
\end{proposition}

\begin{proof}
Let $\mathcal{F}_{0}=\mathcal{R}_{\alpha_{0}},$ $\mathcal{F}_{1}%
=\mathcal{R}_{\alpha_{1}}$ and $g_{1}\left(  k\right)  =k$ for all
$k\in\mathbb{N}.$ Suppose that $g_{n}$ and $\mathcal{F}_{n}$ have been
defined. If $\alpha_{n+1}<\omega$, let $\mathcal{F}_{n+1}=\mathcal{R}%
_{\alpha_{n+1}}$ and $g_{n+1}=g_{n}.$ If $\alpha_{n+1}\geq\omega,$ pick
$x\left(  k,n\right)  \in c_{00}$ for each $k\in\mathbb{N}$ such that

\begin{enumerate}
\item $\min\,$supp $x\left(  k,n\right)  \geq k,$

\item $\left|  \left|  x\left(  k,n\right)  \right|  \right|  _{\ell^{1}%
}=1/\theta_{n+1},$ and

\item $\left|  \left|  x\left(  k,n\right)  \right|  \right|  _{\left[
\mathcal{F}_{n_{1}},\dots,\mathcal{F}_{n_{s}},\mathcal{F}_{0}\right]  }%
\leq\frac{1}{\left|  A\right|  }$ whenever $n_{1},\dots,n_{s}\leq n,$
\end{enumerate}

\noindent where $A = \left\{  \left(  n_{1},\dots,n_{s}\right)  :\theta
_{n_{1}}\dots\theta_{n_{s}}>\theta_{n+1}\right\}  .$ Choose a nondecreasing
function $g_{n+1}:\mathbb{N}\rightarrow\mathbb{N}$ such that $g_{n+1}\geq
g_{n}$ and supp $x\left(  k,p\right)  \subseteq\lbrack k,g_{n+1}\left(
k\right)  )$ for all $1\leq p\leq n,$ $k\in\mathbb{N}.$ Then choose families
$\mathcal{G}_{n+1}$ and $\mathcal{H}_{n+1}$ corresponding to $\alpha_{n+1}$
and $g_{n+1}$ using Lemma \ref{L23}. Finally, define $\mathcal{F}%
_{n+1}=\left(  \mathcal{H}_{n+1},\mathcal{G}_{n+1}\left[  \mathcal{S}%
_{1}^{g_{n+1}}\right]  \right)  .$ Note that $\iota\left(  \mathcal{F}%
_{n}\right)  =\alpha_{n}$ for all $n.$ This completes the inductive definition
of the families $\left(  \mathcal{F}_{n}\right)  _{n=0}^{\infty}.$\newline 

\noindent\textbf{Claim. }If $\alpha_{n+1}\geq\omega,$ then $\left\|  x\left(
k,n\right)  \right\|  \leq2$ for all $k\in\mathbb{N}.$

\noindent Let $x=x\left(  k,n\right)  $ and suppose $\left\|  x\right\|
=\sum_{E\in\mathcal{E}}t\left(  E\right)  \left\|  Ex\right\|  _{\mathcal{F}%
_{0}},$ where $\mathcal{E}$ is the set of all leaves of an $\left(
\mathcal{F}_{n}\right)  $-admissible tree. Take
\[
\mathcal{E}^{\prime}=\left\{  E\in\mathcal{E}:h\left(  E\right)  =\left(
0,n_{1},\dots,n_{s}\right)  ,\text{ }\left(  n_{1},\dots,n_{s}\right)  \in
A\right\}
\]
and $\mathcal{E}^{\prime\prime}=\mathcal{E\setminus E}^{\prime}.$ Now
$E\in\mathcal{E}^{\prime\prime}$ only if $t\left(  E\right)  \leq\theta
_{n+1}.$ Therefore%
\[
\sum_{E\in\mathcal{E}^{\prime\prime}}t\left(  E\right)  \left|  \left|
Ex\right|  \right|  _{\mathcal{F}_{0}}\leq\theta_{n+1}\sum_{E\in
\mathcal{E}^{\prime\prime}}\left|  \left|  Ex\right|  \right|  _{\mathcal{F}%
_{0}}\leq\theta_{n+1}\left\|  x\right\|  _{\ell^{1}}=1.
\]
If $\left(  n_{1},\dots,n_{s}\right)  \in A,$ let $\mathcal{L}_{\left(
n_{1},\dots,n_{s}\right)  }=\left\{  E\in\mathcal{E}^{\prime}:h\left(
E\right)  =\left(  0,n_{1},\dots,n_{s}\right)  \right\}  .$
Now
\[
\sum_{E\in\mathcal{L}_{\left(  n_{1},\dots,n_{s}\right)  }}t\left(  E\right)
\left|  \left|  Ex\right|  \right|  _{\mathcal{F}_{0}}\leq\sum_{E\in
\mathcal{L}_{\left(  n_{1},\dots,n_{s}\right)  }}\left|  \left|  Ex\right|
\right|  _{\mathcal{F}_{0}}\leq\left\|  x\right\|  _{[\mathcal{F}_{n_{1}%
},\dots,\mathcal{F}_{n_{s}},\mathcal{F}_{0}]}\leq\frac{1}{\left|  A\right|  }%
\]
by condition 3. Hence
\[
\sum_{E\in\mathcal{E}^{\prime}}t\left(  E\right)  \left|  \left|  Ex\right|
\right|  _{\mathcal{F}_{0}}\leq\sum_{\left(  n_{1},\dots,n_{s}\right)  \in
A}\frac{1}{\left|  A\right|  }=1.
\]
Thus%
\[
\left\|  x\right\|  =\sum_{E\in\mathcal{E}}t\left(  E\right)  \left\|
Ex\right\|  =\sum_{E\in\mathcal{E}^{\prime}}t\left(  E\right)  \left\|
Ex\right\|  +\sum_{E\in\mathcal{E}^{\prime\prime}}t\left(  E\right)  \left\|
Ex\right\|  \leq2.
\]
This proves the claim.

Since $\alpha_{n}<\sup_{m}\alpha_{m}=\omega^{\omega^{\xi}}$ for all
$n\in\mathbb{N},$ there exist $n_{1}<n_{2}<n_{3}<\dots$ such that $\sup
_{s}\alpha_{n_{s}+1}=\omega^{\omega^{\xi}}$ and $\alpha_{n_{s}+1}\geq\omega$
for all $s\in\mathbb{N}.$ Note that this implies by choice that $\sup_{s}%
\iota\left(  \mathcal{G}_{n_{s}+1}\right)  =\omega^{\omega^{\xi}}.$ Now choose
$q_{1}<q_{2}<q_{3}<\dots$ such that $q_{s+1}>\max$ supp $x\left(  q_{s}%
,n_{r}\right)  ,$ $1\leq r\leq s.$ Let $L=\left\{  q_{1},q_{2},q_{3}%
,\dots\right\}  \in\left[  \mathbb{N}\right]  $ and $q\left(  F\right)
=\left\{  q_{m}:m\in F\right\}  $ for all $F\in\left[  \mathbb{N}\right]
^{<\infty}.$ Define
\[
\mathcal{G}=\left\{  F:s\leq F\text{ and }q\left(  F\right)  \in
\mathcal{G}_{n_{s}+1}\text{ for some }s\in\mathbb{N}\right\}  .
\]
Then $\iota\left(  \mathcal{G}\right)  =\omega^{\omega^{\xi}}.$ For $s\leq m,$
supp $x\left(  q_{m},n_{s}\right)  \subseteq\lbrack q_{m},g_{n_{s}+1}\left(
q_{m}\right)  )\in\mathcal{S}_{1}^{g_{n_{s}+1}}.$ Hence if $s\leq F$,
$q\left(  F\right)  \in\mathcal{G}_{n_{s}+1}$ for some $s\in\mathbb{N},$ and
$x_{m}=\frac{x\left(  q_{m},n_{s}\right)  }{\left\|  x\left(  q_{m}%
,n_{s}\right)  \right\|  }$ for all $m\in F,$ then%
\[
\bigcup_{m\in F}\text{supp}\,x_{m}\in\mathcal{G}_{n_{s}+1}\left[
\mathcal{S}_{1}^{g_{n_{s}+1}}\right]  \subseteq\mathcal{F}_{n_{s}+1}.
\]
Thus, for all $\left(  a_{m}\right)  \in c_{00},$%
\begin{align*}
\bigl\|\sum_{m\in F}a_{m}x_{m}\bigr\|  &  \geq\theta_{n_{s}+1}\bigl\|%
\sum_{m\in F}a_{m}x_{m}\bigr\|_{\mathcal{F}_{n_{s}+1}}\\
&  = \theta_{n_{s}+1}\bigl\|\sum_{m\in F}a_{m}x_{m}\bigr\|_{\mathcal{\ell}%
^{1}}\\
&  \geq\frac{\theta_{n_{s}+1}}{2}\sum_{m\in F}\left|  a_{m}\right|  \left\|
x\left(  q_{m},n_{s}\right)  \right\|  _{\mathcal{\ell}^{1}}\text{ by the
claim,}\\
&  =\frac{1}{2}\sum_{m\in F}|a_{m}| \text{ by condition 2.}%
\end{align*}
Finally, if $F\in\mathcal{F}_{n}\left[  \mathcal{G}\right]  $ for some
$n\in\mathbb{N},$ write $F=\bigcup_{s=1}^{k}F_{s}$ where $F_{1}<\dots<F_{k},$
$F_{s}\in\mathcal{G}$, $1\leq s\leq k,$ and $\left\{  \min F_{1},\dots,\min
F_{k}\right\}  \in\mathcal{F}_{n}.$ For $1\leq s\leq k,$ choose a normalized
sequence $(x_{m})_{m\in F_{s}}$ as above. Now for all $\left(  a_{m}\right)
\in c_{00},$%
\begin{align*}
\left\|  \sum_{m\in F}a_{m}x_{m}\right\|   &  =\left\|  \sum_{j=1}^{k}\left(
\sum_{m\in F_{s}}a_{m}x_{m}\right)  \right\| \\
&  \geq\theta_{n}\sum_{j=1}^{k}\left\|  \sum_{m\in F_{s}}a_{m}x_{m}\right\| \\
&  \geq\frac{\theta_{n}}{2}\sum_{m\in F}\left|  a_{m}\right|  .
\end{align*}
\end{proof}

\section{Standard Schreier families}

For all limit ordinals $\alpha<\omega_{1},$ fix a sequence of ordinals
strictly increasing to $\alpha.$ If $\beta=\omega^{\beta_{1}}\cdot m_{1}%
+\dots+\omega^{\beta_{k}}\cdot m_{k}$ is a limit ordinal, determine
$\mathcal{S}_{\beta}$ using the sequence%
\[
\hat{\beta}_{n}=\left\{
\begin{array}
[c]{cc}%
\omega^{\beta_{1}}\cdot m_{1}+\dots+\omega^{\beta_{k}}\cdot\left(
m_{k}-1\right)  +\omega^{\beta_{k}-1}\cdot n & \text{if }\beta_{k}\text{ is a
successor,}\\
\omega^{\beta_{1}}\cdot m_{1}+\dots+\omega^{\beta_{k}}\cdot\left(
m_{k}-1\right)  +\omega^{\zeta_{n}} & \text{if }\beta_{k}\text{ is a limit,}%
\end{array}
\right.
\]
where $\left(  \zeta_{n}\right)  $ is the chosen sequence of ordinals
increasing to $\beta_{k}.$ It is clear that if $\alpha$ is a countable limit
ordinal such that $\ell\left(  \alpha\right)  \leq\eta$ for some $\eta
<\omega_{1},$ then $\widehat{\left(  \omega^{\eta}\cdot m+\alpha\right)  }%
_{n}=\omega^{\eta}\cdot m+\hat{\alpha}_{n}$ for all $m,n\in\mathbb{N}.$
Throughout this section, we assume that the Schreier families $\mathcal{S}%
_{\alpha}$ are defined using these choices. For such ``standard'' Schreier
families, the converse of Theorem \ref{T2} holds. We begin by establishing
some lemmas.

\begin{lemma}
\label{L25}If $\alpha$ and $\eta$ are countable ordinals such that
$\ell\left(  \alpha\right)  \leq\eta$ and $m \in\mathbb{N}$, then
$\mathcal{S}_{\alpha}\left[  \mathcal{S}_{\omega^{\eta}\cdot m}\right]  =
\mathcal{S}_{\omega^{\eta}\cdot m+\alpha}$.
\end{lemma}

\begin{proof}
The proof is by induction on $\alpha.$ The case $\alpha=0$ is clear. The
result holds for $\alpha=1$ by definition of $\mathcal{S}_{\omega^{\eta}\cdot
m+1}.$ Suppose the lemma is true for some $\alpha.$ Then
\begin{align*}
\mathcal{S}_{\alpha+1}\left[  \mathcal{S}_{\omega^{\eta}\cdot m}\right]   &
=\left(  \mathcal{S}_{1}\left[  \mathcal{S}_{\alpha}\right]  \right)  \left[
\mathcal{S}_{\omega^{\eta}\cdot m}\right]  =\mathcal{S}_{1}\left[
\mathcal{S}_{\alpha}\left[  \mathcal{S}_{\omega^{\eta}\cdot m}\right]  \right]
\\
&  =\mathcal{S}_{1}\left[  \mathcal{S}_{\omega^{\eta}\cdot m+\alpha}\right]
=\mathcal{S}_{\omega^{\eta}\cdot m+\alpha+1}.
\end{align*}
Suppose $\alpha$ is a limit ordinal and the lemma holds for all $\gamma
<\alpha.$ By the remark above, $\omega^{\eta}\cdot m+\hat{\alpha}_{n}%
=\widehat{\left(  \omega^{\eta}\cdot m+\alpha\right)  }_{n}$ for all
$m,n\in\mathbb{N}.$ Now
\begin{align*}
F  &  \in\mathcal{S}_{\alpha}\left[  \mathcal{S}_{\omega^{\eta}\cdot m}\right]
\\
&  \Leftrightarrow F\in\mathcal{S}_{\hat{\alpha}_{n}}\left[  \mathcal{S}%
_{\omega^{\eta}\cdot m}\right]  \text{ for some }n\leq\min F,\\
&  \Leftrightarrow F\in\mathcal{S}_{\omega^{\eta}\cdot m+\hat{\alpha}_{n}%
}\text{ for some }n\leq\min F\text{ by induction,}\\
&  \Leftrightarrow F\in\mathcal{S}_{\omega^{\eta}\cdot m+\alpha}.
\end{align*}
\end{proof}

For the next theorem, fix a countable successor ordinal $\xi$ and a
nondecreasing sequence of ordinals $\left(  \beta_{n}\right)  _{n=1}^{\infty}$
such that $\sup\limits_{n\in\mathbb{N}}\beta_{n}=\omega^{\xi}$ nontrivially.
Also let $\mathcal{F}_{0}$ be a regular family containing $\mathcal{S}_{0}$
such that $\iota\left(  \mathcal{F}_{0}\right)  =\alpha_{0}<\omega
^{\omega^{\xi}},$ and let $\left(  \theta_{n}\right)  _{n=1}^{\infty}$ be a
nonincreasing null sequence in $\left(  0,1\right)  .$ In the present context,
the ordinal $\gamma\left(  \varepsilon,m\right)  $ defined at the beginning of
\S\ref{S4} becomes%
\[
\gamma=\gamma( \varepsilon,m) =\max\{ \ell( \alpha_{0}) +\beta_{n_{s}}%
+\dots+\beta_{n_{1}}: \varepsilon\theta_{n_{1}}\theta_{n_{2}}\dots
\theta_{n_{s}}>\theta_{m}\}
\]
for all $m\in\mathbb{N}$ and $\varepsilon>0$ $(\max\emptyset= 0)$. Denote the
immediate predecessor of $\xi$ by $\xi- 1$.

\begin{theorem}
\label{T4}Follow the notation above and apply the standard choices to define
Schreier families. If there exists $\varepsilon>0$ such that for all
$\beta<\omega^{\xi},$ there exists $m\in\mathbb{N}$ satisfying $\gamma\left(
\varepsilon,m\right)  +2+\beta<\beta_{m},$ then $I_{b}\left(  T\left(
\mathcal{F}_{0},\left(  \theta_{n},\mathcal{S}_{\beta_{n}}\right)
_{n=1}^{\infty}\right)  \right)  =\omega^{\omega^{\xi}\cdot2}$. Otherwise,
$I_{b}\left(  T\left(  \mathcal{F}_{0},\left(  \theta_{n},\mathcal{S}%
_{\beta_{n}}\right)  _{n=1}^{\infty}\right)  \right)  =\omega^{\omega^{\xi}}$.
\end{theorem}

\begin{proof}
If there exists $\varepsilon>0$ with the above properties, then Theorem
\ref{T2} yields that $I_{b}\left(  T\left(  \mathcal{F}_{0},\left(  \theta
_{n},\mathcal{S}_{\beta_{n}}\right)  _{n=1}^{\infty}\right)  \right)
=\omega^{\omega^{\xi}\cdot2}$. Now assume that such $\varepsilon$ does not
exist. Given $\varepsilon>0$, there exists $r=r\left(  \varepsilon\right)
\in\mathbb{N}$ such that for all $m\in\mathbb{N}$, $\gamma\left(
\varepsilon,m\right)  +2+\omega^{\xi-1}\cdot r\geq\beta_{m}.$ Let $m_{0}%
\in\mathbb{N}$ be such that $\beta_{m_{0}}> \ell(\alpha_{0})+2+\omega^{\xi
-1}\cdot r.$ Fix $m\geq m_{0}.$ In particular, $\gamma\left(  \varepsilon
,m\right)  \neq0.$ Hence there exist $n_{1},\dots,n_{s}\in\mathbb{N}$ such
that $\varepsilon\theta_{n_{1}}\dots\theta_{n_{s}}>\theta_{m}$ and
$\ell(\alpha_{0})+\beta_{n_{s}}+\dots+\beta_{n_{1}}+2+\omega^{\xi-1}\cdot
r\geq\beta_{m}.$ Choose $r_{0} \in\mathbb{N}$ such that $\ell(\alpha_{0}) + 2
\leq\omega^{\xi-1}\cdot r_{0}$ and write $\beta_{n}=\omega^{\xi-1}\cdot
r_{n}+\gamma_{n}$ for all $n \in\mathbb{N}$, where $r_{n}\in\mathbb{N}%
\cup\left\{  0\right\}  $ and $\gamma_{n}<\omega^{\xi-1}.$ Then $r_{0}%
+r_{n_{1}}+\dots+r_{n_{s}}+r\geq r_{m}.$ If $r_{n}>0,$%
\begin{align*}
\mathcal{S}_{\beta_{n}}  &  =\mathcal{S}_{\omega^{\xi-1}\cdot r_{n}+\gamma
_{n}}=\mathcal{S}_{\gamma_{n}}\left[  \mathcal{S}_{\omega^{\xi-1}\cdot r_{n}%
}\right]  \text{ by Lemma \ref{L25}}\\
&  \supseteq\mathcal{S}_{\omega^{\xi-1}\cdot r_{n}}.
\end{align*}
The inclusion is obvious if $r_{n} = 0$. Therefore, using Lemma \ref{L25}
again,
\begin{align*}
\left[  \mathcal{S}_{\omega^{\xi-1}\cdot\left(  r_{0}+r+1\right)
},\mathcal{S}_{\beta_{n_{1}}},\dots,\mathcal{S}_{\beta_{n_{s}}}\right]   &
\supseteq\left[  \mathcal{S}_{\omega^{\xi-1}\cdot\left(  r_{0}+ r+1\right)
},\mathcal{S}_{\omega^{\xi-1}\cdot r_{n_{1}}},\dots,\mathcal{S}_{\omega
^{\xi-1}\cdot r_{n_{s}}}\right] \\
&  =\mathcal{S}_{\omega^{\xi-1}\cdot(r_{n_{s}}+\dots+r_{n_{1}}+r_{0}+r+1)}.
\end{align*}
Since $\beta_{m}\leq\omega^{\xi-1}\cdot\left(  r_{0}+r_{n_{1}}+\dots+r_{n_{s}%
}+r+1\right)  ,$ it follows from \cite[Proposition 3.2(a)]{OTW} that there
exists $j_{m}\in\mathbb{N},$ such that
\begin{align*}
\mathcal{S}_{\beta_{m}}\cap\left[  \mathbb{N}_{j_{m}}\right]  ^{<\infty}  &
\subseteq\mathcal{S}_{\omega^{\xi-1}\cdot(r_{n_{s}}+\dots+r_{n_{1}}%
+r_{0}+r+1)}\\
&  \subseteq[ \mathcal{S}_{\omega^{\xi-1}\cdot\left(  r_{0}+r+1\right)
},\mathcal{S}_{\beta_{n_{1}}},\dots,\mathcal{S}_{\beta_{n_{s}}}],
\end{align*}
where $\mathbb{N}_{j}$ is the integer interval $[j,\infty)$ for all
$j\in\mathbb{N}.$
By Proposition \ref{P0}, there exists a sequence $\left(  \ell_{m}\right)
\subseteq\mathbb{N}$ converging to $\infty$ such that, defining $\mathcal{F}%
_{n}$ to be $(\mathcal{S}_{\beta_{n}}\cap[ \mathbb{N}_{\ell_{n}}] ^{<\infty})
\cup\mathcal{S}_{0}$ for all $n \in$, $T\left(  \mathcal{F}_{0},\left(
\theta_{n},\mathcal{S}_{\beta_{n}}\right)  _{n=1}^{\infty}\right)  $ is
isomorphic to $T\left(  \mathcal{F}_{0},\left(  \theta_{n},\mathcal{F}%
_{n}\right)  _{n=1}^{\infty}\right)  $.
Let $k_{m}=\max\{j_{m},\ell_{n_{1}},\dots,\ell_{n_{s}}\}$,
\[
\mathcal{B}_{m}=\left\{  B\in\left[  \mathbb{N}\right]  ^{<\infty}:\ell
_{m}\leq B\text{ and }\left|  B\right|  \leq k_{m}\right\}  ,
\]
and define $\mathcal{H}=( \cup_{m=m_{0}}^{\infty}\mathcal{B}_{m})
\cup\mathcal{S}_{\omega^{\xi-1}\cdot(r_{0}+r+1) }.$ If $m\geq m_{0},$ then
$\mathcal{F}_{m}\subseteq[( \mathcal{H}) ^{2},\mathcal{F}_{n_{1}}%
,\dots,\mathcal{F}_{n_{s}}] .$ Indeed, if $F\in\mathcal{F}_{m},$ then
$F\in\mathcal{S}_{0}$ or $F\in\mathcal{S}_{\beta_{m}}\cap[ \mathbb{N}%
_{\ell_{m}}] ^{<\infty}$. In the former case it is clear that $F\in[(
\mathcal{H}) ^{2},\mathcal{F}_{n_{1}},\dots,\mathcal{F}_{n_{s}}] .$ Suppose
$F\in\mathcal{S}_{\beta_{m}}\cap[ \mathbb{N}_{\ell_{m}}] ^{<\infty}$. Then
$F=F_{1}\cup F_{2},$ where $F_{1} = F\cap\lbrack\ell_{m},k_{m})$ and
$F_{2}=F\setminus F_{1}.$ Clearly $F_{1}\in\mathcal{B}_{m}\subseteq[
\mathcal{H},\mathcal{F}_{n_{1}},\dots,\mathcal{F}_{n_{s}}] $ and
\begin{align*}
F_{2}  &  \in\mathcal{S}_{\beta_{m}}\cap[ \mathbb{N}_{k_{m}}] ^{<\infty}\\
&  \subseteq[ \mathcal{S}_{\omega^{\xi-1}\cdot(r_{0}+r+1) },\mathcal{S}%
_{\beta_{n_{1}}}\cap[ \mathbb{N}_{k_{m}}] ^{<\infty},\dots,\mathcal{S}%
_{\beta_{n_{s}}}\cap[ \mathbb{N}_{k_{m}}] ^{<\infty}]\\
&  \subseteq[ \mathcal{H},\mathcal{F}_{n_{1}},\dots,\mathcal{F}_{n_{s}}] .
\end{align*}
Hence $\mathcal{F}_{m}\subseteq[( \mathcal{H}) ^{2},\mathcal{F}_{n_{1}}%
,\dots,\mathcal{F}_{n_{s}}] .$ This proves that the family $\mathcal{G}%
_{\epsilon}= (\mathcal{H}) ^{2}$ satisfies the hypothesis of Proposition
\ref{P11}. Note that $\iota(( \mathcal{H}) ^{2}) =\iota( \mathcal{H})
\cdot2=\omega^{\omega^{\xi-1}\cdot(r_{0}+r+1) }\cdot2.$ Applying Proposition
\ref{P11}, we obtain
\[
I_{b}\left(  T\left(  \mathcal{F}_{0},\left(  \theta_{n},\mathcal{F}%
_{n}\right)  _{n=1}^{\infty}\right)  \right)  \leq\sup_{\varepsilon>0}%
\sup_{n\in\mathbb{N}}\, [ \omega^{\omega^{\xi-1}\cdot\left(  r_{0}+r\left(
\varepsilon\right)  +1\right)  }\cdot2\cdot\omega^{{\beta_{n}}\cdot\omega}]
=\omega^{\omega^{\xi}}.
\]
Since the reverse inequality holds by Theorem \ref{P7}, the proof is complete.
\end{proof}

It is worthwhile to record the statement of Theorem \ref{T4} for finite
$\beta_{n}$'s.

\begin{corollary}
Suppose that $\mathcal{F}_{0}$ is a regular family containing $\mathcal{S}%
_{0}$ such that $\iota(\mathcal{F}) < \omega^{\omega}$, and that $(\theta
_{n})$ is a nonincreasing null sequence in $(0,1)$ such that $\theta_{n+m}
\geq\theta_{n}\theta_{m}$ for all $n, m \in$. Let $X = T(\mathcal{F}%
_{0},(\theta_{n},\mathcal{S}_{n})^{\infty}_{n=1})$. If $\lim_{m}\limsup
_{n}\theta_{m+n}/\theta_{n} > 0$, then $I(X) = \omega^{\omega\cdot2}$.
Otherwise, $I(X) = \omega^{\omega}$.
\end{corollary}

We conclude by stating without proof a special case of the result when $\xi=
0$. For any $n\in\mathbb{N}$, define $\mathcal{A}_{n}$ to be the family of all
subsets of $\mathbb{N}$ of cardinality $\leq n$.

\begin{proposition}
Suppose that $\mathcal{F}_{0}$ is a regular family containing $\mathcal{S}%
_{0}$ and $\iota(\mathcal{F}_{0}) < \omega$. Let $(k_{n})$ be a sequence in
$\mathbb{N}$ such that $\lim k_{n} = \infty$ and $(\theta_{n})^{\infty}_{n=1}$
be a nonincreasing null sequence in $(0,1)$. Denote the space $T(\mathcal{F}%
_{0},(\theta_{n},\mathcal{A}_{k_{n}}) _{n=1}^{\infty})$ by $Y$. Assume that
every term $(\theta_{n},\mathcal{A}_{k_{n}})$ is essential in the sense that
there exists a nonzero $x \in Y$ such that $\|x\| = \theta_{n}\sum^{k_{n}%
}_{j=1}\|E_{j}x\|$ for some $E_{1} < \dots< E_{k_{n}}$. Then $I_{b}(Y)=\omega$
if
\[
\inf_{r\in\mathbb{N}}\sup\bigl\{\frac{\theta_{m}}{\theta_{n}}: k_{m} \geq
rk_{n}\bigr\} > 0.
\]
Otherwise, $I_{b}(Y)=\omega^{2}$.
\end{proposition}

\end{document}